\documentclass{article}
\title{One- and multi-dimensional CWENOZ reconstructions
	for implementing boundary conditions without ghost cells}
\author{M. Semplice
	\thanks{
		Dipartimento di Scienza e Alta Tecnologia --
		Università dell'Insubria; 
		Via Valleggio, 11 -- 22100 Como (Italy)
		\email{matteo.semplice@uninsubria.it}
	} 
	\and 
	E. Travaglia 
	\thanks{
		Dipartimento di Matematica --
		Università di Torino;
		Via C. Alberto, 10 -- 10124 Torino (Italy)
		\email{elena.travaglia@unito.it}
	}
	\and 
	G. Puppo
	\thanks{
		Dipartimento di Matematica --
		Università La Sapienza;
		P.le Aldo Moro, 5 -- 00185 Roma (Italy)
		\email{gabriella.puppo@uniroma1.it}
	}
}

\newcommand{\keywords}[1]{{\newcommand{\and}{; }\par{\bf Keywords.} #1}}
\newcommand{\subclass}[1]{{\newcommand{\and}{; }\par{\bf MSC.} #1}}
\newenvironment{acknowledgements}{\paragraph{Acknowledgements}}{}
\newcommand{\email}[1]{Email: {\sl #1}}

\usepackage{graphicx}
\graphicspath{{immagini/}{tikzplot/}}
\usepackage{subfig}

\usepackage{tikz}
\usetikzlibrary{decorations.pathreplacing}
\usepackage{pgfplotstable}\usetikzlibrary{shapes,spy,plotmarks}
\usepackage{booktabs}

\usepackage{amsmath,amsfonts}
\usepackage{nicefrac}
\newcommand{\CWENO}[1]{\ensuremath{\mathsf{CWENO#1}}}
\newcommand{\CWENOZ}[1]{\ensuremath{\mathsf{CWENOZ#1}}}
\newcommand{\CWAO}{\ensuremath{\mathsf{CWENO\mbox{-}AO}}}
\newcommand{\CWZAO}{\ensuremath{\mathsf{CWENOZ\mbox{-}AO}}}
\newcommand{\CWb}[1]{\ensuremath{\mathsf{CWb#1}}}
\newcommand{\CWZb}[1]{\ensuremath{\mathsf{CWZb#1}}}

\newcommand{\DX}{\mathrm{\Delta}x}
\newcommand{\DT}{\mathrm{\Delta}t}

\newcommand{\Ogrande}{\mathcal{O}}

\newcommand{\ca}[1]{\overline{#1}}
\newcommand{\Prec}{P_{\text{\sf rec}}}
\newcommand{\Popt}{P_{\text{\sf opt}}}
\newcommand{\dOpt}{d_0} %
\newcommand{\omegaOpt}{\omega_0} %
\newcommand{\OSC}{\mathrm{OSC}}

\begin{document}
\maketitle
	
\begin{abstract}
	We address the issue of point value reconstructions from cell averages
	in the context of third order finite volume schemes, 
	focusing in particular on the cells close to the boundaries of the domain.
	In fact, most techniques known in the literature rely on the creation of ghost cells outside the boundary
	and on some form of extrapolation from the inside that,
	taking into account the boundary conditions, fills the ghost cells with appropriate values,
	so that a standard reconstruction can be applied also in boundary cells.
    In (Naumann, Kolb, Semplice, 2018), motivated by the difficulty of choosing appropriate
	boundary conditions at the internal nodes of a network, a different technique was explored 
	that avoids the use of ghost cells, 
	but instead employs for the boundary cells a different stencil,
	biased towards the interior of the domain.
	
	In this paper, extending that approach,
	which does not make use of ghost cells, 
	we propose a more accurate reconstruction for the one-dimensional case
	and a two-dimensional one for Cartesian grids.
	In several numerical tests we compare the novel reconstruction
	with the standard approach using ghost cells.
\keywords{high order finite volume schemes 
	\and
	boundary conditions without ghost cells
	\and
	hyperbolic systems
	\and
	CWENOZ reconstruction
	\and
	adaptive order reconstructions
}
\subclass{65M08 \and 76M12}
\end{abstract}

\section{Introduction}
\label{sec:intro}

Computing in an efficient way
accurate albeit non-oscillatory solutions of conservation laws
requires the employment of high-order accurate numerical schemes.
Their design encounters the main difficulties in controlling
the spurious oscillations near discontinuities
and near the domain boundaries.
The first is
well tackled by reconstructions of
the Weighted Essentially Non-Oscillatory (WENO) class
\cite{Shu:97} (see the reviews \cite{Shu:2009:WENOreview,Shu:2016:WENOreview,Shu:2020:WENOreview})
or 
Central Weighted Essentially Non-Oscillatory) (CWENO)
\cite{LPR:00:SIAMJSciComp,SCR:cwenoAMR,Baeza:19:CWENOglobalaverageweight,ZhuQiu18:triFV,BGFB:2020} (general results the finite volume setting are proven in
\cite{CPSV:cweno,CSV19:cwenoz,SV:cwao}).

The issue of boundary treatment for hyperbolic conservation laws is usually tackled by constructing ghost points or ghost cells outside the computational domain and by setting their values with appropriate techniques; after this, a high-order non-oscillatory reconstruction procedure can be applied also close to the boundary, despite its large stencil, thanks to the ghost values. This approach is of course delicate, especially with finite-difference discretizations on non-conforming meshes. In this context a very successful technique is the Inverse Lax-Wendroff approach introduced in \cite{TanShu:10:ILW}, rendered more computationally efficient in \cite{TanWWangShuNing:12:SILW}, and further studied and extended for example in
\cite{LiShuZhang:16,LuFangTanShuZhang:16,LuShuTanZhang:21}; a quite up-to-date review may be found in \cite{Shu:17:ILWhandbook}. A modified procedure enhancing its accuracy and stability has been proposed in \cite{ZhaoHuangRuuth:20}.
Other approaches, still based on an inverse Lax-Wendroff procedure but more tailored to coupling conditions on networks can be found in \cite{BorscheKall:14,ContarinoEtAl:16}.
A different approach, entirely based on WENO extrapolation is studied in \cite{BaezaMuletZorio:16jcp,BaezaMuletZorio:16sema}.

In \cite{KolbSemplice:boundary} a different approach was considered. There, in a one-dimensional finite volume context, ghost values are entirely avoided and the point value reconstruction at the boundary is performed with a CWENO type non-oscillatory reconstruction that makes use of only interior cell averages.
The reconstruction stencil for the last cell at the boundary is not symmetric,
but extends only towards the interior of the computational domain.
 Then, the boundary flux is determined from the reconstructed value and the boundary conditions. Achieving non-oscillatory properties when a discontinuity is close to the boundary requires the inclusion of very low degree (down to a constant one, in fact) polynomials in the CWENO procedure; this, in turn, calls for infinitesimal linear weights in order not to degrade the accuracy on smooth solutions. This type of CWENO reconstructions have been studied in general in \cite{SV:cwao}.

In this paper, 
we first enhance the accuracy of the boundary treatment of \cite{KolbSemplice:boundary} 
by employing an Adaptive Order CWENO-Z reconstruction from \cite{SV:cwao} 
and furthermore extend it to two space dimensions.
In particular, in \S\ref{sec:rec1d} 
we describe the new one-dimensional reconstruction that avoids ghost cells
and, in \S\ref{sec:num1d}, compare it with the one of \cite{KolbSemplice:boundary} on numerical tests.
The novel two-dimensional no-ghost reconstruction is then presented in \S\ref{sec:rec2d}
and the corresponding numerical results that compares it with the ghosted approach of \cite{CSV19:cwenoz} are presented in \S\ref{sec:num2d}.
Some final remarks and conclusions are drawn in \S\ref{sec:concl}.

\section{The novel CWENOZb reconstruction in one space dimension}
\label{sec:rec1d}

We start recalling here the operators that define a generic Central WENO reconstruction, which will be useful later.

Central WENO is a procedure to reconstruct point values of a function
from its cell averages; 
it is different from the classical WENO 
by the fact that it performs a single non-linear weight computation per cell
and outputs a polynomial globally defined in the cell,
which is later evaluated at reconstruction points.

In defining a Central WENO, one starts selecting an optimal polynomial, 
denoted here by $\Popt$,
which should be chosen to have the maximal desired accuracy;
the CWENO reconstruction polynomial, in fact, 
will be very close to this one 
when the cell averages in the stencil are a sampling of a smooth enough function.
For the cases when a discontinuity is present in the stencil of $\Popt$,
a sufficient number of alternative polynomials, 
$P_1,\ldots,P_m$,
typically with lower degree and with a smaller stencil, 
should be made available to the blending procedure.

The CWENO operator than computes a nonlinear blending of all polynomials as follows.
First a set of positive {\sl linear} or {\sl optimal coefficients} is chosen,
with the only requirement that
$\dOpt+d_1+\ldots+d_m=1$.
Then, the reconstruction polynomial is defined by
\begin{equation}
\label{eq:CW}
\Prec
=
\omegaOpt\left(\frac{\Popt-\sum_{i=1}^{m}d_iP_i}{\dOpt}\right)
+\sum_{i=1}^{m}\omega_iP_i
\end{equation}
The quantities $\omega_i$ appearing above are called 
{\sl nonlinear weights}; 
when $\omega_i=d_i$ for $i=0,\ldots,n$, then $\Prec=\Popt$ and the
reconstruction will have the maximal accuracy.
When a discontinuity is present in the stencil, the nonlinear weight should deviate from their optimal values in order to avoid the occurrence of spurius oscillations in the numerical scheme.
In practice, the nonlinear weights  are computed with the help of 
oscillation indicators associated to each polynomial, 
that should be $o(1)$ when the polynomial interpolates smooth data and
$\Ogrande(1)$ when the polynomial interpolates discontinuous data.
The construction is independent from the specific form of these indicators,
which here we denote generically as $\OSC[P]$;
typically, the Jiang-Shu indicators from \cite{JiangShu:96} are employed.

When the reconstruction is denoted by 
$\CWENO(\Popt;P_1,\ldots,P_m)$, the nonlinear coefficients are computed as in the
original WENO construction, namely as
\begin{equation} \label{eq:omegas}
    \alpha_k = \frac{d_k}{\left(I_k+\epsilon\right)^p} 
    \qquad
    \omega_k = \frac{\alpha_k}{\sum_j \alpha_j} 
\end{equation}
where $\epsilon$ is a small parameter and $p\geq1$.
For detailed results on the accuracy of such a reconstruction, see
\cite{CPSV:cweno} and the references therein.

Better accuracy on smooth data, especially on coarse grids, 
without sacrificing the non-oscillatory properties,
can be obtained by computing the nonlinear weights as in the
WENOZ construction of \cite{DB:2013}, namely as
\begin{equation} \label{eq:omegaZs}
    \alpha_k = d_k\left[1+ \left(\frac{\tau}{I_k+\epsilon}\right)^p\right] 
    \qquad
    \omega_k = \frac{\alpha_k}{\sum_j \alpha_j}.
\end{equation}
In this case, we denote the reconstruction as
$\CWENOZ(\Popt;P_1,\ldots,P_m)$.
Here above, $\tau$ is quantity that is supposed to 
be much smaller than the individual indicators
when the data in the entire reconstruction stencil are smooth enough.
For efficiency, this {\sl global smoothness indicator} 
should be computed as a linear combination of the other oscillators.
For results on the optimal choices for $\tau$ in a CWENO setting
and the accuracy of the resulting reconstructions, see
\cite{CSV19:cwenoz} and the references therein.

The accuracy results of both \CWENO\ and \CWENOZ\ require that certain
relations among the accuracy of all polynomials involved are satisfied;
the precise conditions for optimal accuracy 
involve also the parameters $\epsilon$ and $p$ \cite{CPSV:cweno,CSV19:cwenoz}, 
but the as a rule of thumb one should always have
$\mathrm{deg}(\Popt)\leq2\mathrm{deg}(P_k)$ for $k=1\ldots,m$.
If controlling spurious oscillations require the inclusion in the nonlinear combination
of polynomials with degree smaller than $\frac12\mathrm{deg}(\Popt)$,
optimal accuracy can still be achieved if the linear weights of these 
additional polynomials of very low degree are infinitesimal, 
i.e. chosen as $\Ogrande(\DX^r)$ for some $r>0$.
In order to distinguish this situation and to easily spot the polynomials with 
infinitesimal linear weights, we adopt for this case the notations
$\CWAO(\Popt;P_1,\ldots,P_m;Q_1,\ldots,Q_m)$,
when \eqref{eq:omegas} is used for the nonlinear weights,
and
$\CWZAO(\Popt;P_1,\ldots,P_m;Q_1,\ldots,Q_m)$,
when \eqref{eq:omegaZs} is used instead.
This was studied on a specific example in \cite{KolbSemplice:boundary} for the CWENO case
and in general for \CWZAO\ in \cite{SV:cwao}.
This latter contains a thorough study of sufficient conditions on $r$ and 
on the other parameters that guarantee optimal convergence rates for a generic \CWZAO\ reconstruction.

\subsection{CWENO-boundary reconstruction of \cite{KolbSemplice:boundary}}
A third-order accurate reconstruction that does not make use of ghost cells has
been introduced in \cite{KolbSemplice:boundary}. The reconstruction coincides
with the \CWENO3 reconstruction of \cite{LPR:00:SIAMJSciComp}, with variable
$\epsilon$ parameter as in \cite{Kolb:14,CS:epsweno,CPSV:cweno}. 
In particular, for the $j$-th cell one considers
the following polynomials:
$P_j^{(2)}$, which is the optimal second degree polynomial interpolating $\ca{u}_{j-1},\ca{u}_j,\ca{u}_{j+1}$, 
$P^{(1)}_{j,L}$ and $P^{(1)}_{j,R}$,
which are the linear polynomials interpolating 
$\ca{u}_{j-1},\ca{u}_j$ and  $\ca{u}_j,\ca{u}_{j+1}$ respectively.		
\CWENO3 is a shorthand for $\CWENO(P_j^{(2)};P^{(1)}_{j,L},P^{(1)}_{j,R})$.
This
reconstruction produces a second degree, uniformly third order accurate,
polynomial defined in each cell, using the cell averages in a stencil of three cells; it can thus be computed on every cell in the domain except for the last one close to each boundary.
	
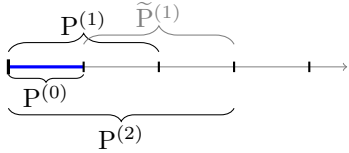
\begin{figure}
	\begin{center}
		\begin{tikzpicture}
\usetikzlibrary{decorations.pathreplacing}

\draw[help lines,->] (0,0) -- (4.5,0);
\draw[very thick,blue] (0,0) --(1,0);

\draw[very thick] (0,-.1) -- (0,.1);
\foreach \x in {1,2,3,4}
  \draw[xshift=\x cm,thick] (0,-.07) -- (0,.07);
 
\draw [yshift=-1mm,decorate,decoration={brace,mirror,amplitude=3pt}]
    (0,0) -- (1,0) node [midway,anchor=north] {$\mathrm{P}^{(0)}$};

\draw [yshift=-5.5mm,decorate,decoration={brace,mirror,amplitude=5pt}]
    (0,0) -- (3,0) node [midway,anchor=north,yshift=-1mm] {$\mathrm{P}^{(2)}$};

\draw [yshift=2mm,decorate,decoration={brace,amplitude=7pt}]
    (0,0) -- (2,0) node [midway,anchor=south,yshift=1mm] {$\mathrm{P}^{(1)}$};

\draw [gray,yshift=3mm,decorate,decoration={brace,amplitude=7pt}]
    (1,0) -- (3,0) node [midway,anchor=south,yshift=1mm] {$\widetilde{\mathrm{P}}^{(1)}$};

\end{tikzpicture}
	\end{center}
	\caption{Illustration of the stencil for the 3-rd order reconstruction in the last cell. Blue: cell where the reconstruction is computed. Black: the polynomials involved in $\Prec$. Gray: additional polynomial for $\tau_{(b3)}$.}
	\label{fig:CWZb3:stencil}
\end{figure}

In the first cell of the domain, the reconstruction is replaced with an	adaptive-order reconstruction $\CWAO(\hat{P}_1^{(2)};P^{(1)_{1,R}};P_1^{(0)})$ in which the stencils of the quadratic $\hat{P}_1^{(2)}$ and of the linear $P^{(1)}$ polynomial do not involve ghost cells (see also Fig.~\ref{fig:CWZb3:stencil}) 
and $P_1^{(0)}$ is the constant polynomial with value $\ca{u}_1$.
In particular $\hat{P}_1^{(2)}$ is the parabola
interpolating the cell averages $\ca{u}_1,\ca{u}_2,\ca{u}_3$. The last cell is treated symmetrically.

The inclusion of the constant polynomial $P^{(0)}$ is necessary to prevent oscillations when a	discontinuity is present one cell away from the boundary and giving it an infinitesimal linear weight allows to guarantee the optimal order of convergence for the reconstruction procedure on smooth data.
More precisely, in \cite{KolbSemplice:boundary} 
it is shown that choosing the linear weights 
as $d^{(0)}=\min(\DX^{\hat{m}},0.01)$ for the constant polynomial, 
$d^{(1)}=0.25$ for the linear one
and consequently setting $d_0=1-d^{(1)}-d^{(0)}$, 
guarantees the optimal accuracy on smooth data when $\hat{m}\in[1,2]$, 
provided $\epsilon=\DX^q$ with $q\geq\hat{m}$.

In general a small $\epsilon$ yields good results on discontinuities, 
but keeping $q=1$ seems desirable to avoid rounding problems in the computation of the nonlinear weights. 
The combination $\hat{m}=2$ and $q=1$, 
however does not fulfill the hypotheses of the convergence result of \cite{KolbSemplice:boundary}; 
in practice, however, the reconstruction appears to give rise 
nevertheless to a third order accurate scheme 
but degraded accuracy can be observed at low grid resolutions. 
As an extreme example in this sense, 
let us consider the linear transport of a periodic initial datum in a periodic domain. 
Of course there would be no need to employ the no-ghost reconstruction in this case, 
since it would be trivial to fill in the ghost values (except maybe for considerations on parallel communication), 
but this example serves quite well to illustrate the situation on smooth data.

In Table~\ref{tab:cwb3errors} we report the 1-norm errors observed for the transport of $u(x,0)=\sin(\pi x - \sin(\pi x)/\pi)$	after one period (for full details on the numerical scheme, the reader is referred to the beginning of \S\ref{sec:num1d}). It is evident that for the reconstruction of \cite{KolbSemplice:boundary}, third order error rates are observed only on very fine grids when $d^{(0)}\sim\DX$; for $d^{(0)}\sim\DX^2$, the optimal rate predicted by the theory is observed in practice, but the errors are still larger than its ghosted \CWENO3 counterpart.

\begin{table}
\caption{Errors on the linear transport of $\sin(\pi x - \sin(\pi x)/\pi)$ in a periodic domain, using \CWENO{3} and \CWb{3} reconstructions ($\epsilon=\DX^2$).
}
\label{tab:cwb3errors}
\centering
\footnotesize

\pgfplotstabletypeset[
	col sep=space,
	sci zerofill,
	empty cells with={--},
	every head row/.style=
 {before row=\toprule
		& \multicolumn{2}{c|}{\CWENO{3}} 
		& \multicolumn{2}{c|}{\CWb{3}, $d^{(0)}=\DX$} 
		& \multicolumn{2}{c|}{\CWb{3}, $d^{(0)}=\DX^2$} 
		\\,
 after row=\midrule
 },
	every last row/.style={after row=\bottomrule},
	columns={N,1,r1,5,r5,9,r9},
	create on use/N/.style=
 {create col/expr={2/\thisrow{0}}},
	columns/N/.style={column name={N},precision=0,column type=r|},
	columns/1/.style={column name={error}},
	create on use/r1/.style=
 {create col/dyadic refinement rate={1}},
	columns/r1/.style=
      {column name={rate},fixed zerofill,column type=r|},
	columns/5/.style={column name={error}},
	create on use/r5/.style=
 {create col/dyadic refinement rate={5}},
	columns/r5/.style=
 {column name={rate},fixed zerofill,column type=r|},
	columns/9/.style={column name={error}},
	create on use/r9/.style=
 {create col/dyadic refinement rate={9}},
	columns/r9/.style=
 {column name={rate},fixed zerofill,column type=r|},
	]
	{dati1d/lintraSinCW.err}
\end{table}	
	
\subsection{The novel \CWZb3 reconstruction}
The loss of accuracy at low resolution can be traced back to the relative inability of the smoothness indicators alone to detect a smooth flow on coarse grids.
The net effect is that, when the grid is coarse,
the nonlinear weight of the constant polynomial
in the first and last cells is larger than it would be strictly needed,
degrading the accuracy of the reconstruction there; the errors are then transported into the domain by the flow.

This issue can be successfully counteracted, even on coarse grids, by the employment of Z-weights in the construction. In fact, we recall that the idea behind WENO-Z, see \cite{DB:2013},
is to replace the standard WENO nonlinear weight computation
\eqref{eq:omegas} with \eqref{eq:omegaZs}
where the global smoothness indicator $\tau$ is supposed to be $\tau=o(I_k)$ if the cell averages represent a locally smooth data in the stencil. The improved performances of WENO-Z over WENO, and of \CWENOZ\ over \CWENO\ reconstructions are in fact linked to the superior ability of detecting smooth transitions, already at low grid resolution, granted by the global smoothness indicator $\tau$.
Moreover, the ability of detecting a smooth flow even at low grid resolution
depends on how small is $\tau$ on smooth data;
thus the goal in the optimal design of $\tau$ 
is to choose the coefficients of the linear combination 
$\lambda_0\OSC[\Popt]+\sum_{i=1}^{n}\lambda_k\OSC[P_k] = \Ogrande(\DX^s)$
that maximize $s$
when the data in the stencil of the reconstruction is a sampling of a smooth function \cite{CSV19:cwenoz}.

Our proposal thus consists in defining the new \CWZb3\ reconstruction to coincide with 
$\CWENOZ3=\CWENOZ(P_j^{(2)};P^{(1)}_{j,L},P^{(1)}_{j,R})$ in the domain interior, 
with the adaptive-order reconstruction 
$\CWZAO(\hat{P}_1^{(2)};P_{1,R}^{(1)};P_1^{(0)})$ in the first cell
and with
$\CWZAO(\hat{P}_N^{(2)};P_{N,L}^{(1)};P_N^{(0)})$ in the last cell

We now have to specify our choice of $\tau$.	
The definition of the Jiang-Shu 
oscillation indicators \cite{JiangShu:96} is 
\[
\OSC[P] 
:= 
\sum_{\ell\geq1}
\DX^{2\ell-1}
\int_{\Omega_j}
\left(\frac{\mathrm{d}^{\ell}}{\mathrm{d}x^{\ell}} P \right)^2
\mathrm{dx}	
\]
where $\Omega_j$ is the cell where the reconstruction is applied.

On smooth data, in the domain interior, we have that
\begin{align*}
\OSC[P_j^{(2)}] &= (u'(x_j))^2 \DX^2 + \Ogrande(\DX^4)
\\
\OSC[P^{(1)}_{j,L}] &= (u'(x_j))^2 \DX^2 - u'(x_j)u''(x_j)\DX^3  +\Ogrande(\DX^4)
\\
\OSC[P^{(1)}_{j,R}] &= (u'(x_j))^2 \DX^2 + u'(x_j)u''(x_j)\DX^3 +\Ogrande(\DX^4).
\end{align*}	
so that the combination
\begin{equation}
\label{eq:tau}
\tau_j=\left|2\OSC[P_j^{(2)}]-\OSC[P^{(1)}_{j,L}]-\OSC[P^{(1)}_{j,R}]\right|
\end{equation}
is $\Ogrande(\DX^4)$;
this very low $\tau$ biases very strongly the non-linear weights \eqref{eq:omegaZs}
towards the optimal ones whenever the flow is smooth.
In \cite{CSV19:cwenoz} it is shown that this is the optimal choice and that
it is not possible to obtain a combination of the indicators that is
$o(\DX^4)$ in the third order setup.

We now need to specify a suitable $\tau_1$ for the uncentered stencil of the first cell and $\tau_N$ for the last one. 
Recall that the role of $\tau$ is to indicate whether the data are smooth in the stencil, 
which is composed by the first three cells adjacent to the boundary.
As argued in \cite{SV:cwao}, 
only the polynomials with degree at least one are useful in the construction of $\tau$.
Here we could use the oscillators
of the parabola $P^{(2)}$ fitting the three cell averages $\ca{u}_1,\ca{u}_2,\ca{u}_3$ and the linear polynomial $P^{(1)}_1$ interpolating the first two,
\begin{align*}
\OSC[P^{(2)}] &= (u'(x_j))^2 \DX^2 + \Ogrande(\DX^4)
\\
\OSC[P^{(1)}] &= (u'(x_j))^2 \DX^2 + u'(x_j)u''(x_j)\DX^3  +\Ogrande(\DX^4)
\end{align*}
and thus we cannot exploit the symmetry 
to obtain a global smoothness indicator of size $\Ogrande(\DX^4)$.

Using $\tau$ of $\Ogrande(\DX^3)$ however could make the reconstruction in the boundary cell less performing than the one in the domain interior.
In order to overcome this difficulty, 
one could employ, in the construction of $\tau$, also the indicator of the linear polynomial $\widetilde{P}^{(1)}$ interpolating the averages $\ca{u}_2,\ca{u}_3$.
However, since the role of $\tau$ is to detect smooth flows in the global stencil, 
which is composed by the first three cells,
that is the same reconstruction stencil employed by the second cell, 
a simpler solution (which also allows to save some computations)
is to take instead for the first cell 
the same value of $\tau$ that was computed in the second cell; 
this is $\Ogrande(\DX^4)$ on smooth flows and yields a better reconstruction.

The novel reconstruction procedure that we propose is thus:
\begin{itemize}
\item in all cells except the first and last one, compute the \CWENOZ3 reconstruction polynomial with the optimal definition \eqref{eq:tau} 
of $\tau_j$, 
as in \cite{CSV19:cwenoz};
\item in the first cell, 
apply $\CWZAO(\hat{P}_1^{(2)};P_{1,R}^{(1)};P_1^{(0)})$
with $\tau_1:=\tau_2$
\item in the last cell,
apply
$\CWZAO(\hat{P}_N^{(2)};P_{N,L}^{(1)};P_N^{(0)})$ 
with $\tau_N:=\tau_{N-1}$.
\end{itemize}
After the analysis of \S3.1.1 of \cite{SV:cwao},
it is expected that this reconstruction has the third order of accuracy
for $d^{(0)}=\Ogrande(\DX)$
provided that $p\geq1$ and $\epsilon=\Ogrande(\DX^{\hat{m}})$ for $\hat{m}\in[1,3]$.

As discussed in \cite{CSV19:cwenoz},
the choice of parameters within the allowed ranges
can trade better accuracy on smooth flows 
(larger $\hat{m}$ or smaller $q$)
with a smaller spurious oscillations on discontinuities
(smaller $\hat{m}$ or larger $q$).
In \cite{CSV19:cwenoz} it was found that a good overall choice for \CWENOZ3 was $p=1$ and $\hat{m}=2$ 
and we will adopt these values in all our numerical tests. 
Regarding the infinitesimal linear weight,
the choices $d^{(0)}=\DX$ and $d^{(0)}=\DX^2$
will be compared.

\section{One-dimensional numerical tests}
\label{sec:num1d}
All tests in this section are conducted with a finite volume scheme 
constructed with the method of lines, 
the Local Lax-Friedrichs numerical flux, 
and the third order TVD-RK3.
The \CWENO3 and the \CWENOZ{3} reconstruction 
in the first/last cell make use of one ghost cell outside each boundary, 
which is filled according to the boundary conditions
before computing the reconstruction.
In the same cells, the \CWb3 and the \CWZb3 reconstructions, instead, 
do not make use of ghost cells 
but extend their stencil for one extra cells inwards with respect to their ghosted counterparts.

In both cases, the flux on the boundary face is computed by applying a consistent numerical flux 
(here the Local Lax Friedrichs)
to an inner value determined by the reconstruction 
and an outer value determined by the boundary conditions.
More precisely, 
for periodic boundary conditions, 
the outer value on the left is copied from the inner value at the right boundary and viceversa;
for reflecting boundary conditions in gasdynamics, 
the outer value is the same as the inner one
but has the opposite sign for the velocity;
for Dirichlet boundary, 
the outer value is set to the exact value of the boundary function at time $t_n+c_i\DT$ 
for the $i$-th stage of the Runge-Kutta scheme.

The CFL number is set to $0.45$ in all tests. The numerical tests have been performed with the open-source code {\tt claw1dArena}, see \cite{claw1dArena}.

\subsection{Linear transport}
\paragraph{Periodic solution}
We consider again the linear transport equation $u_t+u_x=0$ in the domain $[-1,1]$ with periodic boundary conditions. We evolve for one period the initial data $u_0(x)=\sin(\pi x - \sin(\pi x)/\pi)$, which has a critical point of order 1 (see \cite{HAP:2005:mappedWENO}), with the \CWENOZ{3} and \CWZb{3} reconstructions.

\begin{table}
\caption{Errors on the linear transport of $\sin(\pi x - \sin(\pi x)/\pi)$ in a periodic domain, using \CWENOZ{3} and \CWZb{3} reconstructions.
	}
\label{tab:cwzb3errors}
\centering
\footnotesize
\pgfplotstabletypeset[
	col sep=space,
	sci zerofill,
	empty cells with={--},
	every head row/.style=
	{before row=\toprule
		& \multicolumn{2}{c|}{\CWENOZ{3}} 
		& \multicolumn{2}{c|}{\CWZb{3}, $d^{(0)}=\DX$} 
		& \multicolumn{2}{c|}{\CWZb{3}, $d^{(0)}=\DX^2$} 
		\\,
		after row=\midrule
	},
	every last row/.style={after row=\bottomrule},
	columns={N,1,r1,5,r5,9,r9},
	create on use/N/.style=
	{create col/expr={2/\thisrow{0}}},
	columns/N/.style={column name={N},precision=0,column type=r|},
	columns/1/.style={column name={error}},
	create on use/r1/.style=
	{create col/dyadic refinement rate={1}},
	columns/r1/.style=
	{column name={rate},fixed zerofill,column type=r|},
	columns/5/.style={column name={error}},
	create on use/r5/.style=
	{create col/dyadic refinement rate={5}},
	columns/r5/.style=
	{column name={rate},fixed zerofill,column type=r|},
	columns/9/.style={column name={error}},
	create on use/r9/.style=
	{create col/dyadic refinement rate={9}},
	columns/r9/.style=
	{column name={rate},fixed zerofill,column type=r|},
	]
	{dati1d/lintraSinCWZ.err}	
\end{table}	

Table~\ref{tab:cwzb3errors} shows that,
according to the results of \cite{SV:cwao},
\CWZb{3} can reach the optimal convergence rate already with $d^{(0)}\sim\DX$ and that the errors obtained without using ghosts are very close to those of the ghosted reconstruction \CWENOZ{3}. As already pointed out in \cite{CSV19:cwenoz}, also here we observe that using Z-weights in \CWENO\ yields lower errors compared to the companion reconstructions with Jiang-Shu weights (compare Tab.~\ref{tab:cwb3errors}).

\paragraph{Smooth solution with time-dependent Dirichlet data}
For this second test, we consider again the linear transport equation on the domain $[-1,1]$, 
but this time we apply time-dependent Dirichlet boundary data on the left (inflow) imposing $u(-1,t)=0.25 - 0.5\sin(\pi(1.0+t)$
and free-flow conditions on the (outflow) boundary at $x=1$.
We start with $u_0(x)=0.25 + 0.5\sin(\pi x)$ and compare the computed cell averages with the exact solution $u(t,x)=u_0(x-t)$.
The final time is set to $1$.
This test was proposed in \cite{TanShu:10:ILW}.

\begin{table}
	\caption{Errors on the smooth linear transport test
		with time-dependent Dirichlet data.
	}
	\label{tab:smootherrors}
	\centering
	\footnotesize
	\pgfplotstabletypeset[
	col sep=space,
	sci zerofill,
	empty cells with={--},
	every head row/.style=
	{before row=\toprule
		& \multicolumn{2}{c|}{\CWZb{3}, $\delta=\DX$} 
		& \multicolumn{2}{c|}{\CWZb{3}, $\delta=\DX^2$} 
		\\,
		after row=\midrule
	},
	every last row/.style={after row=\bottomrule},
	columns={N,5,r5,9,r9},
	create on use/N/.style=
	{create col/expr={2/\thisrow{0}}},
	columns/N/.style={column name={N},precision=0,column type=r|},
	columns/5/.style={column name={error}},
	create on use/r5/.style=
	{create col/dyadic refinement rate={5}},
	columns/r5/.style=
	{column name={rate},fixed zerofill,column type=r|},
	columns/9/.style={column name={error}},
	create on use/r9/.style=
	{create col/dyadic refinement rate={9}},
	columns/r9/.style=
	{column name={rate},fixed zerofill,column type=r|},
	]
	{dati1d/lintraTanShuSmoothCWZ.err}
\end{table}	

The results reported in Tab.~\ref{tab:smootherrors} show that the \CWZb3 reconstruction yields third order error rates already on coarse grids and with $d^{(0)}\sim\DX$. No advantage is seen for the choice $d^{(0)}\sim\DX^2$.

We point out that applying a reconstruction that makes use of ghost cells, 
like \CWENO3 or \CWENOZ3,
would not be straightforward in this case.
In \cite{CGAD:95:RKboundary} 
was observed that accuracy would be 
capped at second order if 
the ghost cell values for the $i$-th stage were to be set
by reflecting the inner ones in 
the exact boundary data at time $t_n+c_i\DT$,
where $c_i$ denotes the abscissa of the $i$-th stage 
of the Runge-Kutta scheme.
In the same paper, 
also a suitable modification of the boundary data that preserve the accuracy of the Runge-Kutta scheme is also proposed.
On the other hand, 
we point out that with the \CWb3 and \CWZb3 reconstructions
this issue of filling the ghost cells is not present
and that the exact boundary data can be employed in the numerical flux computation
without observing losses of accuracy.

\paragraph{Discontinous solution}
Next we consider the same setup of the previous test, 
but impose the boundary value 
\[
u(-1,t)=
\begin{cases}
0.25 &, t\leq1\\
-1 &, t>1
\end{cases}
,
\]
thus introducing a jump in the exact solution at $t=1$. This test was proposed in \cite{TanShu:10:ILW} and, as there, we compute the flow until $t=1.5$.

\begin{figure}
\begin{center}
	\includegraphics[height=4.3cm]{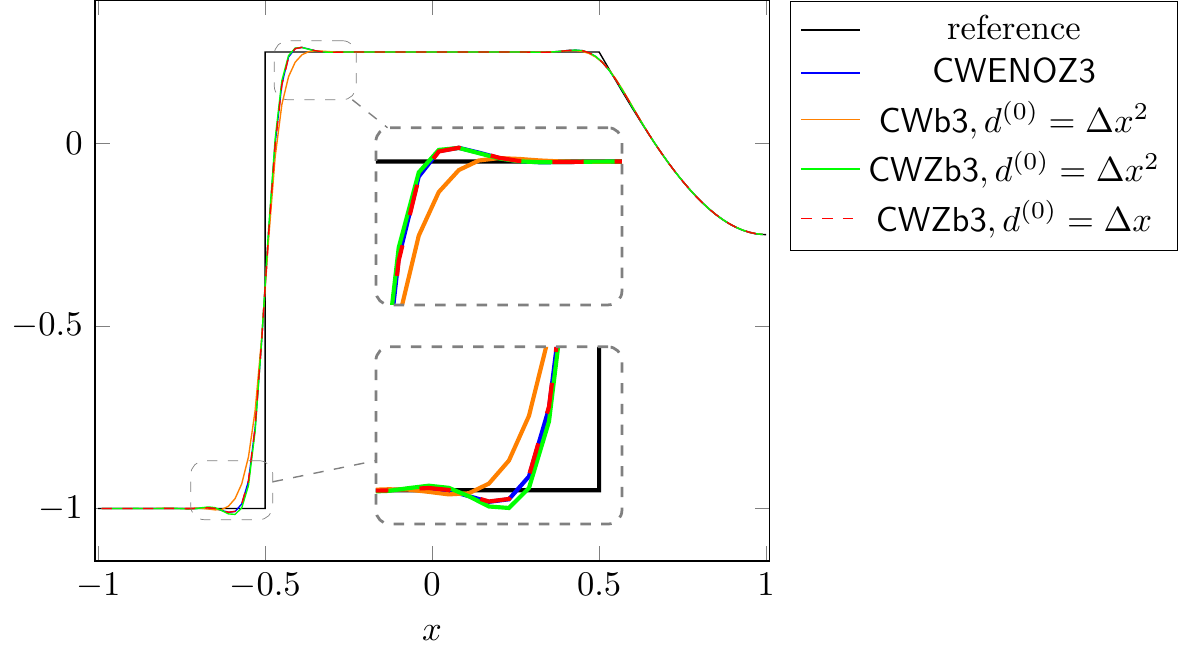}
\end{center}
\caption{Solutions computed with $100$ cells for the discontinuous linear transport test, using \CWENOZ{3}, \CWb3 and \CWZb{3} reconstructions ($\epsilon=\DX^2$).}
\label{fig:tanshudisc}
\end{figure}

The computed solutions are shown in Fig.~\ref{fig:tanshudisc}, 
where we compare the solution computed with  \CWENOZ3 using ghosts and
the no-ghost \CWb3 and \CWZb3.
No difference can be seen in the corner point at $x=0.5$, 
which is originated by a continous but not differentiable boundary data.
On the other hand, the jump at $x=-0.5$ in the final solution
is generated by the discontinuity in the boundary data.
The numerical solution around this jump has
slightly more pronounced oscillations when using \CWZb3 and $d^{(0)}=\DX^2$
and a more smoothed profile when using \CWb3;
\CWZb3 with $d^{(0)}=\DX$ produces an almost idential solution to the one computed by the ghosted \CWENOZ3 reconstruction.

\subsection{Burgers' equation}
For a nonlinear scalar test, 
we consider the Burgers' equation $u_t+(u^2)_x=0$ 
with initial data $u_0(x)=1-\sin(\pi x)$ with periodic boundary conditions, 
so that a shock forms, travels to the right and is located exactly on the boundary at $t=1$.

\begin{figure}	
\begin{center}
		\includegraphics[height=4.3cm]{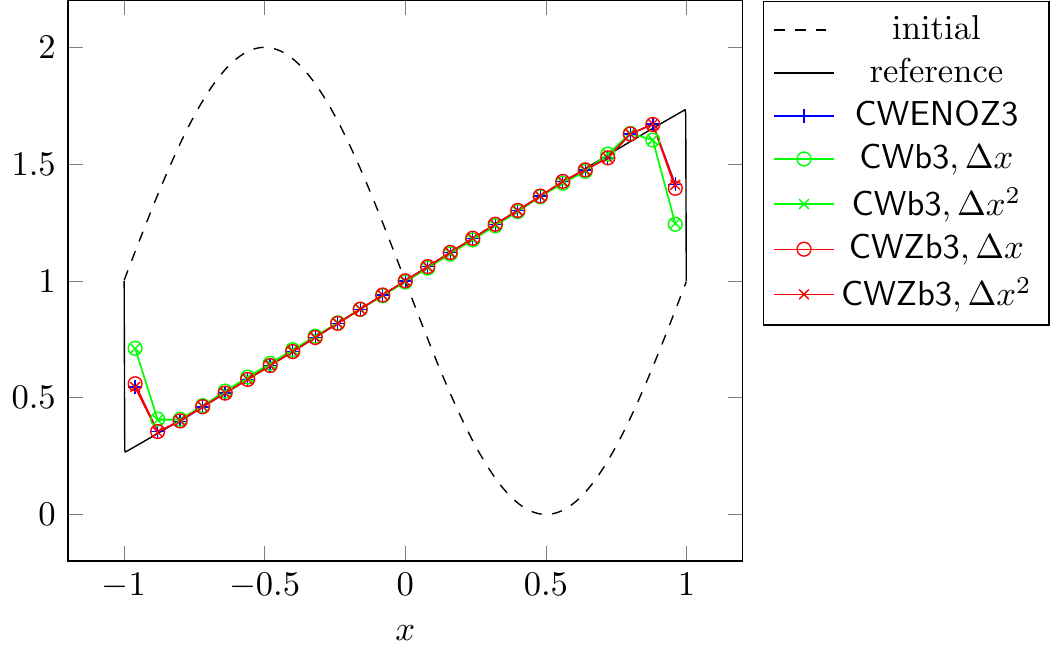}
\end{center}
\caption{Burgers' test with $25$ cells at $t=1$.}
	\label{fig:burgers:3}
\end{figure}

In Fig. \ref{fig:burgers:3} we compare the solutions computed with $25$ cells.
One can see that \CWZb3\ computes a solution which is almost exactly superimposed on the \CWENOZ3, 
despite the fact that using the correct periodic ghost values 
should be an advantage in this test.
The \CWb3 solution is slightly more diffusive
and, in both the no-ghost reconstructions, 
both choices of $d^{(0)}$ yield similar solutions.

\subsection{Euler gas dynamics}
\paragraph{Incoming wave from the left}
In this test we consider a gas initially at rest, 
with $\rho=1, p=1, u=0$ everywhere.
Through a time-dependent Dirichlet boundary condition on the left, 
we introduce the following disturbance 
\[
  \rho(t,0)=1.0+\delta(t)
  \quad
  p(t,0)=1.0+\gamma\delta(t)
  \quad
  \delta(t)=
  \begin{cases}
	0.01(\sin(2\pi t))^3 &, t\in[0,0.5]\\
	0 &, t>0.5
  \end{cases}
\]
The boundary introduces a smooth wave travelling right.
Wall boundary conditions are imposed on the right
and the final time is set at $t=1.25$,
when the wave is being reflected back from the wall.

\begin{figure}
\begin{center}
\includegraphics[height=4.3cm]{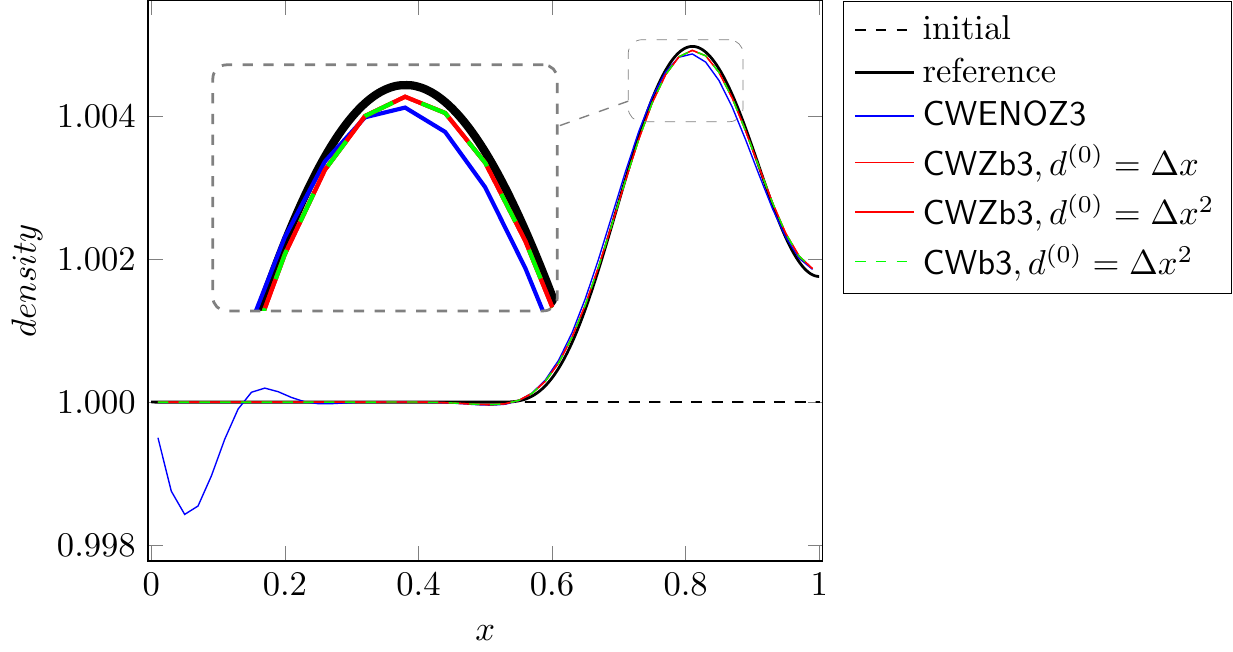}
\end{center}
	\caption{Gasdynamics. Incoming wave test using 50 cells. 
		The reference is computed with 2000 cells and the minmod slope limiter.}
	\label{fig:haysam:3}
\end{figure}

In Fig. \ref{fig:haysam:3}
we report the solutions at time $t=1.25$
computed on 50 cells with the third order ghosted and ghost-free reconstructions,
together with a reference solution computed on 10000 cells with a second order TVD scheme.

Spurious oscillations coming from the Dirichlet boundary conditions on the left side 
are completely absent when using \CWb3 or \CWZb3\ instead of \CWENOZ3. 
Also, a slightly better resolution is observed near the top of the wave.
Here again, we stress that the \CWb3 and the \CWZb3 solutions
have been computed by entirely neglecting the boundary conditions in the reconstruction phase
and passing the exact Dirichlet value at $t_n+c_i\DT$ to the numerical flux 
as outer data on the left boundary.

\paragraph{Sod's Riemann problem}
In this test we use the initial data of the Sod problem, but
we impose wall boundary conditions on both sides.

\begin{figure}
\begin{center}
	\includegraphics[height=4.3cm]{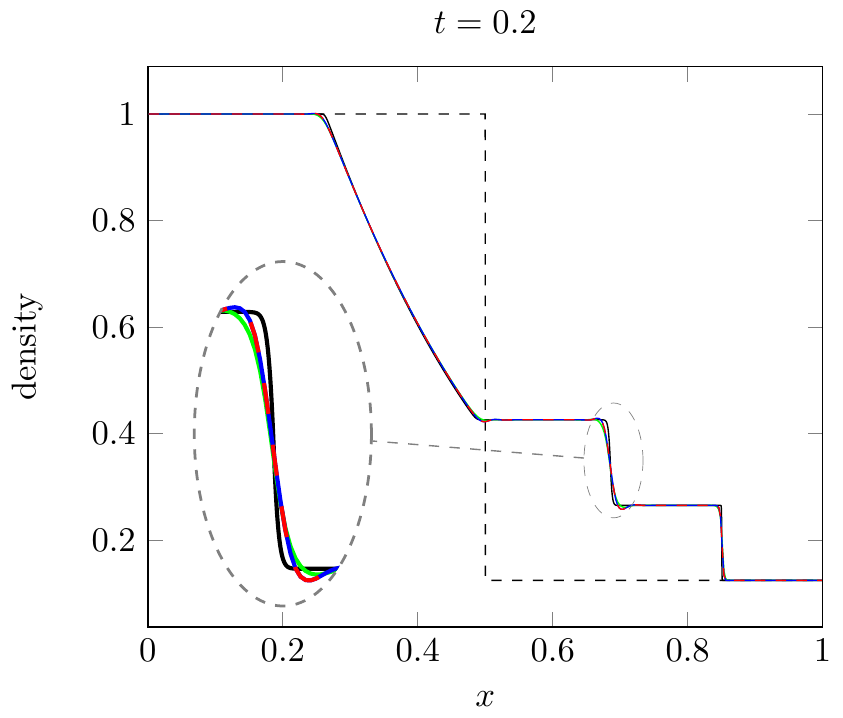}
	\hfill
	\includegraphics[height=4.3cm]{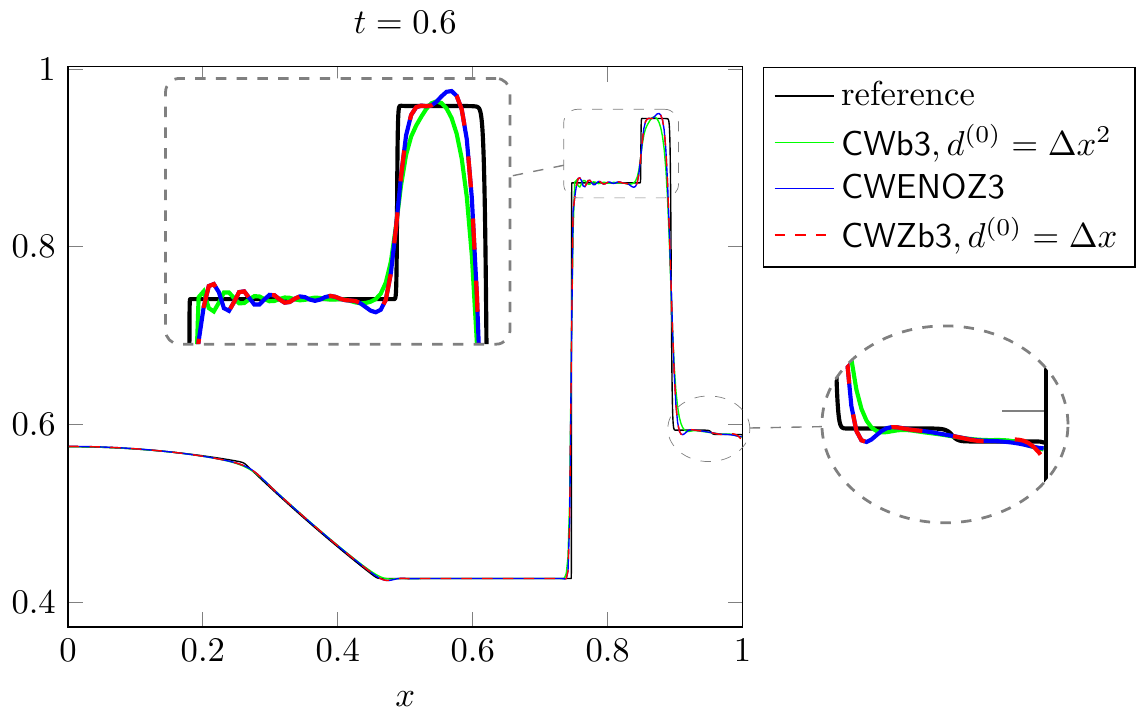}
\end{center}
\caption{Sod test with wall boundary conditions on 400 cells.
	The reference is computed with 10000 cells and a linear reconstruction with the minmod limiter.
}
\label{fig:sodwall1:3}
\end{figure}

In Fig.~\ref{fig:sodwall1:3} we show, in the left panel, the solution at time $t=0.2$, which is before the waves reach the wall; 
the expected solution is thus the usual one.
All three solutions are very close to each other and only a slight extra diffusion can be noticed for the reconstruction that is using the Jiang-Shu nonlinear weights instead of the Z-weights.

Letting the flow evolve past $t=0.2$, 
the shock impinges on the wall and bounces back,
interacting with the right-moving contact around $t=0.41$;
this in turn generates a left-moving shock, 
a very slow contact
and a quite weak right-moving shock;
the right-moving shock then bounces back from the wall 
and interacts with the contact at around $t=0.56$, 
giving rise again to another 
shock-contact-shock interaction pattern.
 
In the right panel of Fig.~\ref{fig:sodwall1:3}
we show the solution at $t=0.6$, and, counting from the left,
we see a rarefaction which is reflecting in the left wall,
two left-moving shocks,
a very slow contact (with speed $0.001$)
and a very weak right-moving shock (density jump below $0.01$).
It can be appreciated that \CWZb3, without using ghost cells, 
computes almost the same solution as the ghosted \CWENOZ3.
As in other tests, \CWb3 is more diffusive.
The very weak shock is barely captured at this resolution
and even the reference solution almost misses it.
The oscillations in the plateaux between the two left-moving shocks
and the hump left of the contact
could be controlled with local characteristic projections,
which was not employed in this computation.

\begin{figure}
\begin{center}
	\includegraphics[height=4.3cm]{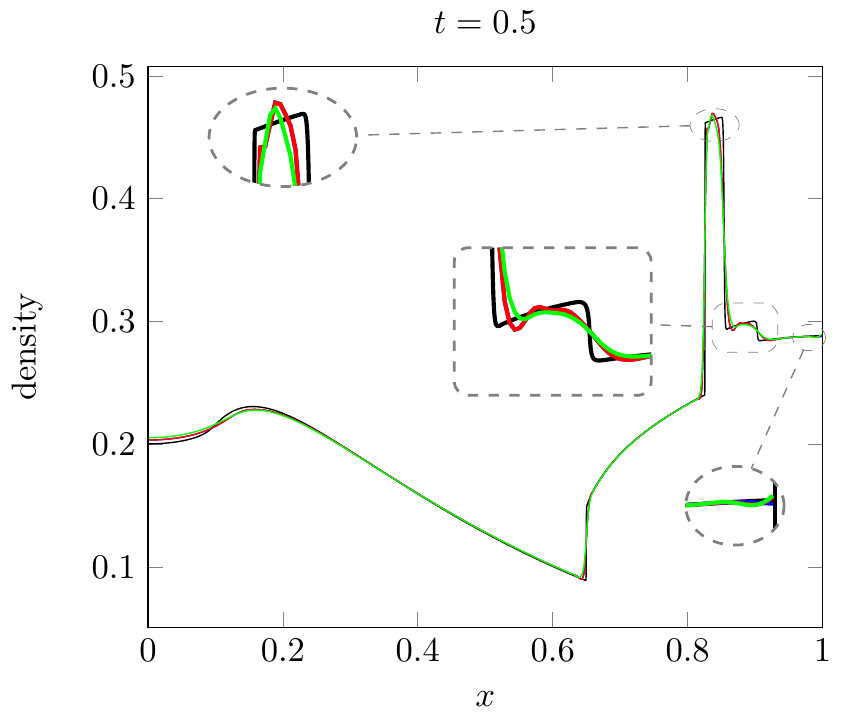}
	\hfill
	\includegraphics[height=4.3cm]{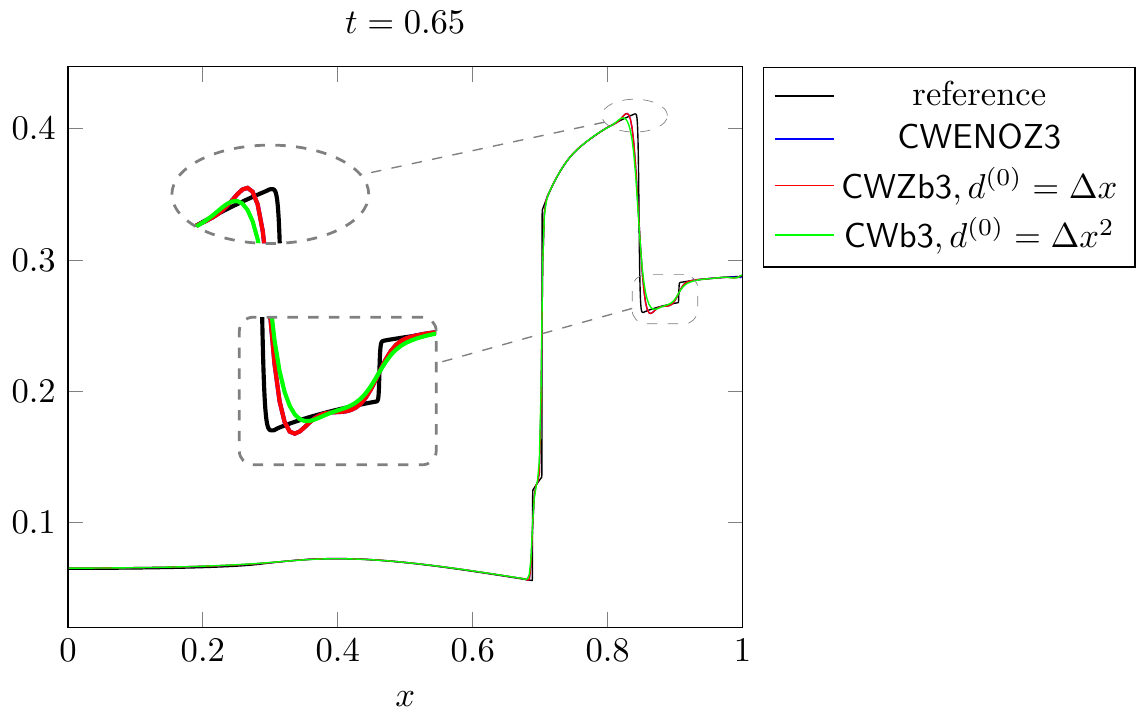}
\end{center}
\caption{Sod test in 3D, with 400 cells. The reference is computed with 4000 cells and a linear reconstruction with the minmod limiter.}
\label{fig:sodwall3:3}
\end{figure}

Finally, we consider the $d$-dimensional version of the same problem. 
Following \cite{Toro:book}, in spherical symmetry this amounts to adding to the Euler equations 
the source term $S(\rho,u,p)=-\tfrac{d-1}{x}[\rho u,\rho u^2, up]^{\mathsf{T}}$.
In particular we show in Fig.~\ref{fig:sodwall3:3}
the solution for $d=3$ at $t=0.5$ and at $t=0.65$.
In this test, 
the source term contribution is computed in each cell
with a two-point gaussian quadrature,
which is fed by the reconstructed values.
We thus test the CWENO-based reconstructions' capability
of easily computing reconstructed values inside the cells.
For all waves, we observe again that \CWZb3 and  \CWENOZ3
produce very similar solutions,
with \CWb3 being slightly more diffusive.

\section{Two-dimensional scheme}
\label{sec:rec2d}

\begin{figure}
\begin{center}
\begin{tikzpicture}
\usetikzlibrary{patterns}
\definecolor{stencil}{gray}{0.8}
\colorlet{cella}{red}

\begin{scope}
\filldraw[fill=stencil](-1.5,-1.5) rectangle (1.5,1.5);
\fill[pattern=north east lines, pattern color=cella] (-.5,-.5) rectangle (.5,.5);
\foreach \x in {-1.5,-0.5,0.5,1.5}
{
  \draw (\x,-1.6) -- (\x,1.6);
  \draw (-1.6,\x) -- (1.6,\x);
}
\draw [blue, xshift=0.15cm,yshift=0.15cm] 
(0,0) circle (1pt) 
-- (0,1) circle (1pt)
-- (1,1) circle (1pt)
-- (1,0) circle (1pt)
-- (0,0);
\draw [blue, xshift=-0.15cm,yshift=0.15cm] 
(0,0) circle (1pt) 
-- (0,1) circle (1pt)
 -- (-1,1) circle (1pt)
 -- (-1,0) circle (1pt)
-- (0,0);
\draw [blue, xshift=0.15cm,yshift=-0.15cm] 
(0,0) circle (1pt) 
-- (0,-1) circle (1pt)
-- (1,-1) circle (1pt)
-- (1,0) circle (1pt)
-- (0,0);
\draw [blue, xshift=-0.15cm,yshift=-0.15cm] 
(0,0) circle (1pt) 
-- (0,-1) circle (1pt)
 -- (-1,-1) circle (1pt)
 -- (-1,0) circle (1pt)
-- (0,0);
\end{scope}

\begin{scope}[xshift=4cm,yshift=-1cm]
\filldraw[fill=stencil](-1.5,-0.5) rectangle (1.5,2.5);
\fill[pattern=north east lines, pattern color=cella] (-.5,-.5) rectangle (.5,.5);
\foreach \x in {-1.5,-0.5,0.5,1.5}
  \draw (\x,-.5) -- (\x,2.6);
\foreach \y in {-0.5,0.5,1.5,2.5}
  \draw (-1.6,\y) -- (1.6,\y);
\draw [very thick] (-1.6,-0.5) -- (1.6,-0.5);
\draw [blue, xshift=0.15cm,yshift=0.15cm] 
(0,0) circle (1pt) 
-- (0,1) circle (1pt)
-- (1,1) circle (1pt)
-- (1,0) circle (1pt)
-- (0,0);
\draw [blue, xshift=-0.15cm,yshift=0.15cm] 
(0,0) circle (1pt) 
-- (0,1) circle (1pt)
 -- (-1,1) circle (1pt)
 -- (-1,0) circle (1pt)
-- (0,0);
\begin{scope}[blue,xshift=0.15cm,yshift=-0.15cm] 
  \draw [dashed]  (0,0) circle (1pt) -- (1,0) circle (1pt);
  \draw                   (0,0) circle (1pt)    (1,0) circle (1pt);
\end{scope}
\begin{scope}[blue,xshift=-0.15cm,yshift=-0.15cm] 
  \draw [dashed]  (0,0) circle (1pt) -- (-1,0) circle (1pt);
  \draw                   (0,0) circle (1pt)    (-1,0) circle (1pt);
\end{scope}
\end{scope}

\begin{scope}[xshift=7cm,yshift=-1cm]
\filldraw[fill=stencil](-0.5,-0.5) rectangle (2.5,2.5);
\fill[pattern=north east lines, pattern color=cella] (-.5,-.5) rectangle (.5,.5);
\foreach \x in {-0.5,-0.5,0.5,1.5,2.5}
  \draw (\x,-.5) -- (\x,2.6);
\draw[very thick] (-.5,-.5) -- (-.5,2.6);
\foreach \y in {-0.5,0.5,1.5,2.5}
  \draw (-0.5,\y) -- (2.6,\y);
\draw[very thick] (-0.5,-.5) -- (2.6,-.5);
\draw [blue, xshift=0.15cm,yshift=0.15cm] 
(0,0) circle (1pt) 
-- (0,1) circle (1pt)
-- (1,1) circle (1pt)
-- (1,0) circle (1pt)
-- (0,0);
\draw [blue, dashed,xshift=0.15cm,yshift=-0.15cm] 
(0,0) circle (1pt) 
-- (1,0) circle (1pt);
\draw [blue,xshift=0.15cm,yshift=-0.15cm] 
(0,0) circle (1pt) 
(1,0) circle (1pt);
\begin{scope}[blue,xshift=0.15cm,yshift=-0.15cm] 
  \draw [dashed]  (0,0) circle (1pt) -- (1,0) circle (1pt);
  \draw                   (0,0) circle (1pt)    (1,0) circle (1pt);
\end{scope}
\begin{scope}[blue,xshift=-0.15cm,yshift=0.15cm] 
  \draw [dashed]  (0,0) circle (1pt) -- (0,1) circle (1pt);
  \draw                   (0,0) circle (1pt)    (0,1) circle (1pt);
\end{scope}
\fill[blue,xshift=-0.15cm,yshift=-0.15cm] (0,0) circle (1.5pt);
\end{scope}
\end{tikzpicture} 
\end{center}
\caption{Stencils for the 2d reconstruction in the middle of the domain (left), 
	at a domain edge (center)
	and at a domain corner (right).
}
\label{fig:stencil:2d}
\end{figure}
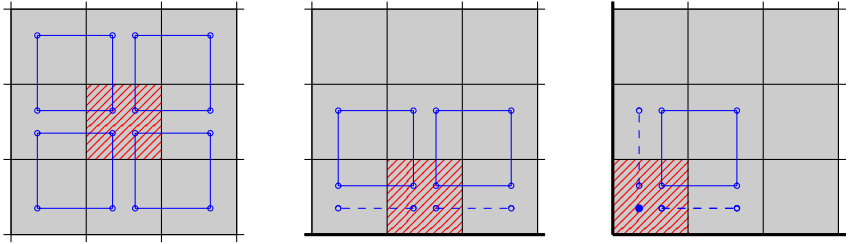

In this section we consider a two-dimensional Cartesian grid, with cells of size $\DX$.
We denote the cells as $\Omega_{i,j}$, 
with the pair of integers $(i,j)$ referring to their position in the grid.
As in the one dimensional case, 
the solution is advanced in time with 
the third order TVD-SSP Runge-Kutta scheme;
the numerical fluxes on each edge of a cell 
are obtained with the two-point gaussian quadrature, 
with point values computed with a two-point numerical flux
fed with the reconstructed values on each side of the edge.
On boundaries, only the inner point value is computed from the reconstruction,
while the outer one is computed according to the boundary conditions.
For example, on a solid wall boundary, the outer value is equal to the inner one, 
except for the normal velocity, which is given the opposite sign.

The reconstruction from cell averages to point values in two space dimensions 
is not obtained by dimensional splitting, 
but is computed by blending polynomials in two spatial variables with a \CWENOZ\ or a \CWZAO\ construction.
The reconstruction operator is called only once per cell
and the polynomial returned is then evaluated at the eight reconstruction points
where the numerical fluxes have to be computed.

Let $\Omega_{i,j}$ be the cell in which the reconstruction is being computed.
In every cell, the reconstruction is computed by a \CWENOZ\ operator 
with optimal polynomial $\Popt$ of degree 2 in two spatial variables (6 degrees of freedom)
associated with a $3\times3$ stencil containing $\Omega_{i,j}$
(see later for the definition of 
the polynomial associated to a stencil).
The reconstruction stencils are depicted in Fig.~\ref{fig:stencil:2d}.
In all panels, the cell in which the reconstruction is being computed is hatched, 
while the stencil of the optimal polynomial of degree 2 is shaded.

The \CWENOZ\ operator is fully specified after the low degree polynomials 
and the global smoothness indicator $\tau$
are also chosen.
The stencils of the low degree polynomials are indicated by circles
joined by solid or dashed lines in Fig.~\ref{fig:stencil:2d}.

In the bulk of the computational domain, 
the reconstruction coincides with the two-dimensional \CWENOZ3 described in \cite{CSV19:cwenoz};
it is defined as a nonlinear combination of second and first degree polynomials:
\[
	\CWENOZ(\Popt;P_{NE},P_{SE},P_{SW},P_{NW}).
\]
The optimal polynomial is associated with the $3\times3$ stencil of cells centered at $\Omega_{i,j}$ (left panel in Fig.~\ref{fig:stencil:2d}).
The four polynomials $P_{NE},P_{SE},P_{SW}$ and $P_{NW}$
are linear polynomials in two variables associated to the four stencils
depicted with solid lines in the figure.
For example, $P_{NE}$ is associate to the stencil 
composed by the cells $\Omega_{r,s}$ 
for $r\in\{i,i+1\}$ and $s\in\{j,j+1\}$.
As in \cite{CSV19:cwenoz}, we define the global smoothness indicator by
\[
	\tau = \big|
	4\OSC[\Popt] - \OSC[P_{NE}] -\OSC[P_{SE}] - \OSC[P_{SW}] - \OSC[P_{NW}]
	\big|,
\]
where $\OSC[P]$ is the multidimensional Jiang-Shu smoothness indicator, 
as defined in \cite{HuShu:WENOtri}.
The nonlinear weights are computed by \eqref{eq:omegaZs} 
starting from the linear weights 
$\dOpt=\nicefrac34$ and 
$d_{NE}=d_{SE}=d_{SW}=d_{NW}=\nicefrac{1}{16}$.

Next we consider the case of a cell adjacent to a domain boundary.
We focus in particular on the case of the bottom boundary, 
which is depicted in the central panel of Fig.~\ref{fig:stencil:2d}.
Here the reconstruction is
$$\CWENOZ(\Popt;P_{NE},P_{NW},\tilde{P}_{E},\tilde{P}_{W}),$$
where $P_{NE}$ and $P_{NW}$ are defined as in the domain bulk.
The stencil of $\Popt$ is biased towards the interior of the domain
and is composed by the cells $\Omega_{r,s}$ 
for $r\in\{i-1,i,i+1\}$ and $s\in\{j,j+1,j+2\}$.
The other two polynomials, $\tilde{P}_{E}$ and $\tilde{P}_{W}$ are degree 1 polynomials
that depend only on the tangential variable, $x$ in the example, and that are constant in 
the direction normal to the boundary. Their stencils are indicated with dashed lines in the figure.
The global smoothness indicator $\tau$ for the cell in the example is copied from the cell $\Omega_{i,j+1}$.
The linear weights are similar to the bulk case, i.e.
$\dOpt=\nicefrac34$ and 
$d_{NE}=d_{NW}=d_{E}=d_{W}=\nicefrac{1}{16}$.
The case of the other boundaries is obtained from this one by symmetry.

Finally we describe the reconstruction in a domain corner,
focusing on the case of the south-west one,
which is represented in the right panel of  Fig.~\ref{fig:stencil:2d}.
Here, for stability purposes, we must include also a constant polynomial
in the \CWENOZ\ operator, denoted with $\tilde{P}_0$,in order to avoid spurious oscillations when a strong wave hits the corner. $\tilde{P}_0$ has of course the constant value coinciding
with the cell average of the corner cell and its 1-cell stencil is represented by the
filled circle in the picture.
Following \cite{SV:cwao}, we assign to the constant polynomial 
and to the $\tilde{P}_{E},\tilde{P}_{N}$ polynomials
an infinitesimal
weight of $d_0=d_N=d_E=\DX^{2}$ 
and the reconstruction in the south-west corner cell is
$$\CWZAO(\Popt;_{NE};\tilde{P}_{E},\tilde{P}_{N},\tilde{P}_0).$$
The stencil of $\Popt$ is again biased towards the interior of the domain
and is composed by the cells $\Omega_{r,s}$ 
for $r\in\{i,i+1,i+2\}$ and $s\in\{j,j+1,j+2\}$.
$\tilde{P}_{E}$, similarly to the previous case, is a degree 1 polynomial
that is constant in the $y$ direction,
while $\tilde{P}_{N}$ is a degree 1 polynomial
that is constant in the $x$ direction.
The global smoothness indicator $\tau$ for the cell in the example is copied from the cell $\Omega_{i+1,j+1}$.
The case of the other corners is obtained from this one by symmetry.

The polynomials associated to the stencils are computed as follows.
Let $\mathcal{S}$ be a collection of neighbours of the cell $\Omega_{i,j}$
that includes the cell itself
and let $\mathrm{\Pi}\subset\mathbb{P}^d(x,y)$
be the subspace of the polynomials of degree $d$ in two spatial variables
where $P_{\mathcal{S}}$ is sought.
If the stencil $\mathcal{S}$ contains 
as many cells as the dimension of $\mathrm{\Pi}$,
the polynomial $P_{\mathcal{S}}$
is the solution of the linear system composed by the equations
$\langle P_{\mathcal{S}} \rangle_{r,s}=\ca{u}_{r,s}$
for all $(r,s)\in\mathcal{S}$,
where the operator $\langle\cdot\rangle_{r,s}$ denotes
the cell average of its argument over the cell $\Omega_{r,s}$.
In the examples above, all polynomials 
with a tilde in their name are computed in this way.

When the cardinality of $\mathcal{S}$ is larger than $\mathrm{dim}(\mathrm{\Pi})$,
we associate to $\mathcal{S}$ the solution of the following constrained least-squares problem:
\begin{equation} \label{eq:ps}
P_{\mathcal{S}}
=
\arg \min 
\left\{ 
\sum_{(r,s) \in \mathcal{S}} 
\big\vert \langle P_{\mathcal{S}} \rangle_{r,s} - \ca{u}_{r,s} \big\vert^2, \, 
\text{such that }  
P_{\mathcal{S}}\in\mathrm{\Pi},
\langle  P_{\mathcal{S}} \rangle_{i,j}= \ca{u}_{i,j} \right \}.
\end{equation}
In the examples above, the polynomials
$\Popt,P_{NE},P_{SE},P_{SW},P_{NW}$
are computed in this way.

On  Cartesian grids,
the constrained least square problem can be easily 
turned into an unconstrained one
by choosing a basis of $\Pi$ 
consisting of a constant function
and of polynomials orthogonal to the constant one.
Explicit expressions for the coefficients
of the polynomials in the domain interior 
can be found in \cite{CWENOandaluz}.

\section{Two-dimensional tests}
\label{sec:num2d}

The numerical scheme has been implemented with the help of the PETSc libraries \cite{petsc-efficient,petsc-user-ref}
for grid management and parallel communications;
the tests were run on a multi-core desktop machine 
equipped with an Intel Core i7-9700 processor
and 64Gb of RAM.
We show the results obtained with the Local Lax-Friedrichs numerical flux.

In all the tests we consider the two-dimensional Euler equations of gas dynamics:
\[ \partial_t \left( \begin{array}{c}
\rho \\ \rho u \\ \rho v \\ E
\end{array}\right) +
\partial_x \left( \begin{array}{c}
\rho u \\ \rho u^2 + p \\ \rho u v \\ u(E+p)
\end{array}\right) +
\partial_y \left( \begin{array}{c}
\rho v \\ \rho u v \\ \rho v^2 + p \\ v(E+p)
\end{array}\right)= 0,
\]
where $\rho$, $u$, $v$, $p$ and $E$ are the density, velocity in $x$ and $y$ direction, pressure and energy per unit mass. 
We consider the perfect gas equation of state
$
E = \frac{p}{\gamma-1} + \frac12 \rho (u^2+v^2)
$
with $\gamma=1.4$.

\subsection{Convergence test}
We compare the novel reconstruction with the one of \cite{CSV19:cwenoz} that makes use of ghost cells on the isentropic vortex test \cite{Shu:97}.
Of course there would be no need to use a ghost-less reconstruction with periodic boundary conditions, since it would be trivial to set up and fill in the ghost cells, but we conduct this as a stress-test to verify the order of accuracy of the novel reconstruction.

The initial condition is characterized by a uniform ambient flow with constant temperature, density, velocity and pressure 
$T_\infty=\rho_\infty=u_\infty=v_\infty=p_\infty = 1.0$, 
onto which the following isentropic perturbations are added in velocity and temperature: 
$$
(\delta u, \delta v) = {\frac {\beta} {2 \pi}} \exp \left( {\frac {1-r^2} {2}} \right) (-y,x), \quad
\delta T = - { \frac {(\gamma - 1 ) \beta^2} {8 \gamma \pi^2}} \exp \left( {1-r^2} \right),
$$
where $r=\sqrt{x^2+y^2}$
and the strength of the vortex is set to $\beta=5.0$.
The computational domain is the square $[-5,5]^2$ with periodic boundary conditions and the final time is set to $t=10$
so that the final exact solution is the same as the initial state.

\begin{table}
	\caption{Errors on the isentropic vortex test,
		 using \CWENOZ{3} and \CWZb{3} reconstructions.
	}
	\label{tab:cwzb3:2derrors}
	\centering
	\footnotesize
	\begin{filecontents*}{dati2d/conv/LLFboth.csv}
50 3.2794e-01 7.7642e-01 7.9875e-01 1.8306e+00 3.0402e-01 7.0536e-01 7.3611e-01 1.7112e+00
100 6.4050e-02 1.3872e-01 1.3985e-01 3.0787e-01 6.2054e-02 1.3306e-01 1.3412e-01 2.9730e-01
200 9.0340e-03 1.9399e-02 1.9509e-02 4.2367e-02 8.8905e-03 1.9026e-02 1.9083e-02 4.1696e-02
400 1.1498e-03 2.4776e-03 2.4899e-03 5.3886e-03 1.1439e-03 2.4620e-03 2.4661e-03 5.3676e-03
800 1.4410e-04 3.1532e-04 3.1688e-04 6.8239e-04 1.4426e-04 3.1535e-04 3.1571e-04 6.8386e-04
1600 1.8019e-05 4.2469e-05 4.2663e-05 9.1225e-05 1.8082e-05 4.2557e-05 4.2614e-05 9.1545e-05

\end{filecontents*}
\pgfplotstabletypeset[
	col sep=space,
	sci zerofill,
	empty cells with={--},
	every head row/.style=
	{before row=\toprule %
		&& \multicolumn{5}{c|}{\CWENOZ{3}} 
		& \multicolumn{6}{c|}{\CWZb{3}} 
		\\,
		after row=\midrule
	},
	every last row/.style={after row=\bottomrule},
	columns={0,1,r1,4,r4,5,r5,8,r8},
	columns/0/.style={column name={N},precision=0,column type=r|},
	create on use/DX/.style=
	{create col/expr={10/\thisrow{0}}},
	columns/DX/.style={column name={$\DX$},sci, precision=2,column type=c|},
	columns/1/.style={column name={density},sci,dec sep align},
	create on use/r1/.style=
	{create col/dyadic refinement rate={1}},
	columns/r1/.style=
	{column name={rate},fixed zerofill,column type=r|},
	columns/4/.style={column name={energy},sci,dec sep align},
	create on use/r4/.style=
	{create col/dyadic refinement rate={4}},
	columns/r4/.style=
	{column name={rate},fixed zerofill,column type=r|},
	columns/5/.style={column name={density},sci,dec sep align},
	create on use/r5/.style=
	{create col/dyadic refinement rate={5}},
	columns/r5/.style=
	{column name={rate},fixed zerofill,column type=r|},
	columns/8/.style={column name={energy},sci,dec sep align},
	create on use/r8/.style=
	{create col/dyadic refinement rate={8}},
	columns/r8/.style=
	{column name={rate},fixed zerofill,column type=r|},
	]
	{dati2d/conv/LLFboth.csv}	
\end{table}	

	We observe third order convergence rates in all variables 
	(1-norm errors in density and energy are shown in Tab.~\ref{tab:cwzb3:2derrors}).
	Compared with the \CWENOZ3 scheme, the errors are no worse, and in some cases slightly better.
	
\subsection{Two-dimensional Riemann problem}
We have run a number of Riemann problems,
in particular configurations B, G and K from \cite{Euler:RP2d:class},
in order to compare the performances of the novel reconstruction
on flows with waves (almost) orthogonal to the boundary.

In our numerical experiments we have noticed that 
choosing correctly the linear weights for the low degree polynomials in the boundary cells
is important to avoid spurious waves and features generated by an anomalous diffusion in the tangential direction; 
this latter shows up for example 
when choosing infinitesimal weights for the planes 
with two cells in the stencil in the middle panel of Fig.~\ref{fig:stencil:2d}.
Since it is on contact waves that spurious diffusion can accumulate over time,
we report only a comparison of the solutions 
computed with the ghosted and the no-ghost reconstruction 
on configuration B of \cite{Euler:RP2d:class},
which involves four contact discontinuities.

\begin{figure}
\begin{center}
	\includegraphics[width=\linewidth]{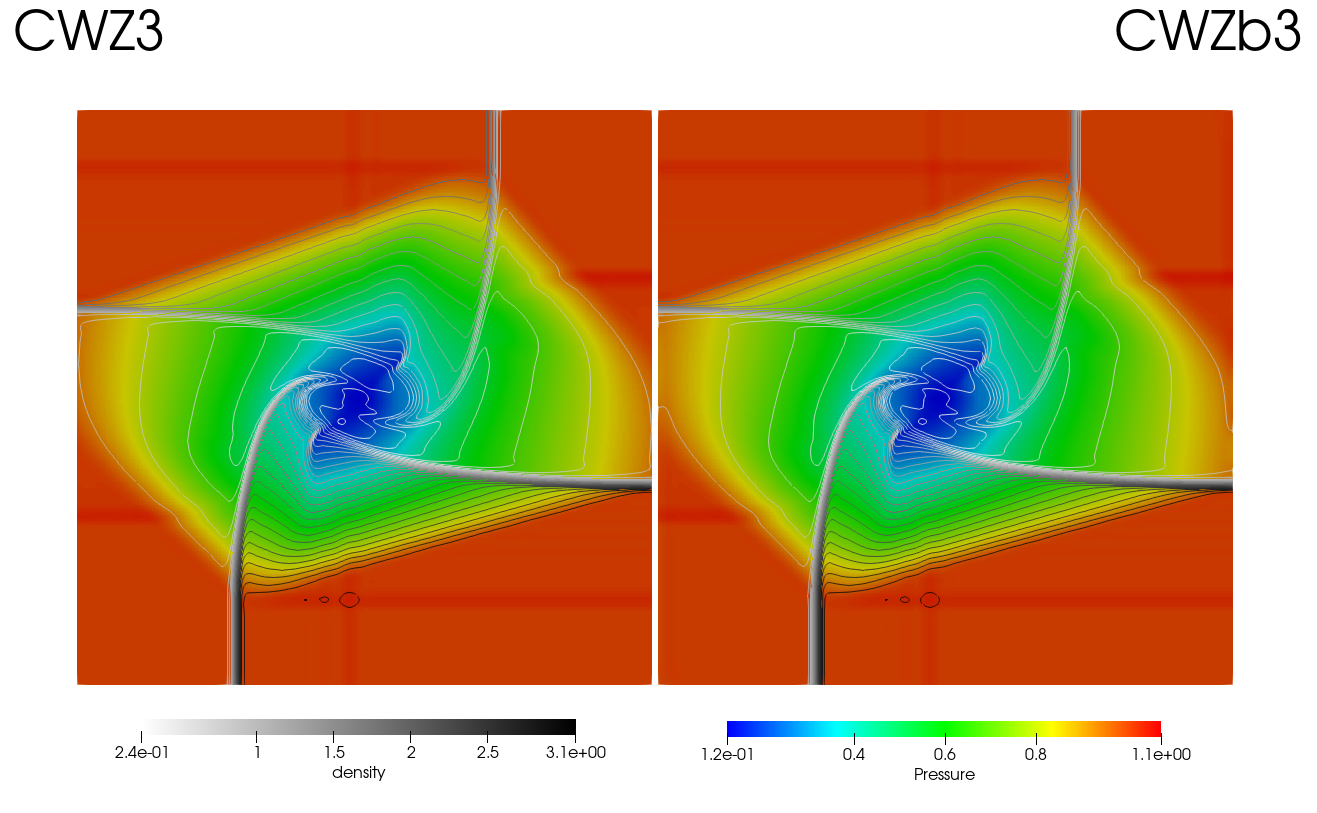}
\end{center}
\caption{Two-dimensional Riemann problem with four slip lines, computed with (left) and without (right) ghost cells.
The colorbar is for pressure;
there are 29 contour lines for density, spaced by $0.1$, from $0.25$ (center) to $3.05$ (in the bottom right part).}
\label{fig:rpb}
\end{figure}

We evolved, in the domain $[-0.5,0.5]^2$ with free-flow boundary conditions, 
an initial configuration with constant data in the four quadrants;
in particular, we set $p=1$ everywhere and
\[
(\rho,u,v)=
\left\{
\begin{array}{c|c}
\text{upper left} & \text{upper right} \\
(2.0,0.75,0.5) & (1.0,0.75,-0.5) \\
\hline
(1.0,-0.75,0.5) & (3.0,-0.75,-0.5) \\
\text{lower left} & \text{lower right} \\
\end{array}
\right\}
\]
so that the solution contains four contact waves rotating in the clock-wise direction.

The solutions computed with and without ghost cells are shown in Fig.~\ref{fig:rpb}.
In the plot the colors stand for pressure (rainbow colorbar) 
and we also showing contour lines of density (grayscale colorbar).
We are focusing on contact waves as they are a good indicator of numerical diffusion, 
since on this kind of waves its effects accumulate over time.
No difference is visible between the two computed solutions, 
indicating that the reconstruction that does not make use of ghost cells
does not introduce significant differences with respect to the standard approach
that makes use of ghosts.
In particular, no wave deformation is visible close to the boundary,
indicating that, with our choice of linear weights,
no extra tangential diffusion is introduced in the boundary cells 
with respect to the cells that are located more inwards into the domain.

\subsection{Radial Sod test}	
Next we run the cylindrical Sod shock tube problem in two space dimensions.
The initial conditions for velocity component are $u,v=0$ everywhere, while density 
and pressure are $(\rho_H,p_H)=(1,1)$ for the central region, i.e. where $r=\sqrt{x^2+y^2}<0.5$, and $(\rho_L,p_L)=(0.125,0.1)$ elsewhere. 
The computational domain is set to $\Omega=[0;1]^2$
with symmetry boundary conditions along $x=0$ and $y=0$ 
and wall boundary conditions on $x=1$ and $y=1$.

\begin{figure}
\begin{center}
	\includegraphics[width=\linewidth]{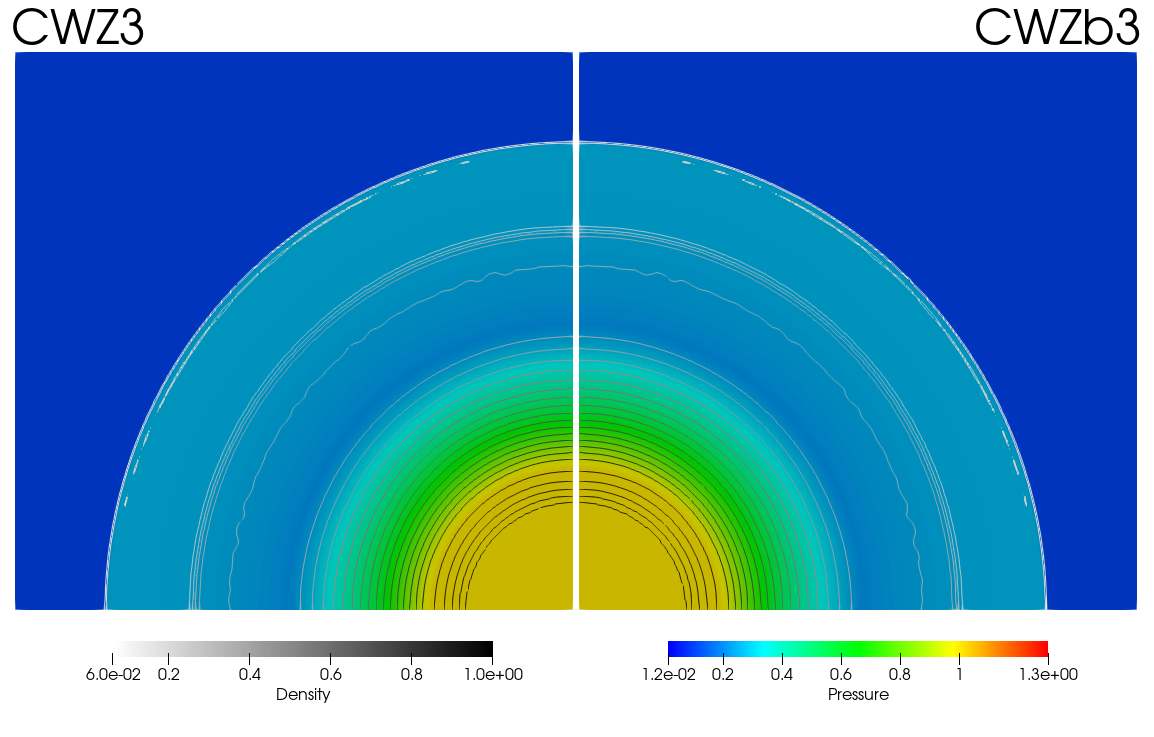}
\end{center}
\caption{Radial Sod solutions at $t=0.2$}
\label{fig:sod:2d}
\end{figure}

\begin{figure}
\begin{center}
  \includegraphics[width=0.32\linewidth]{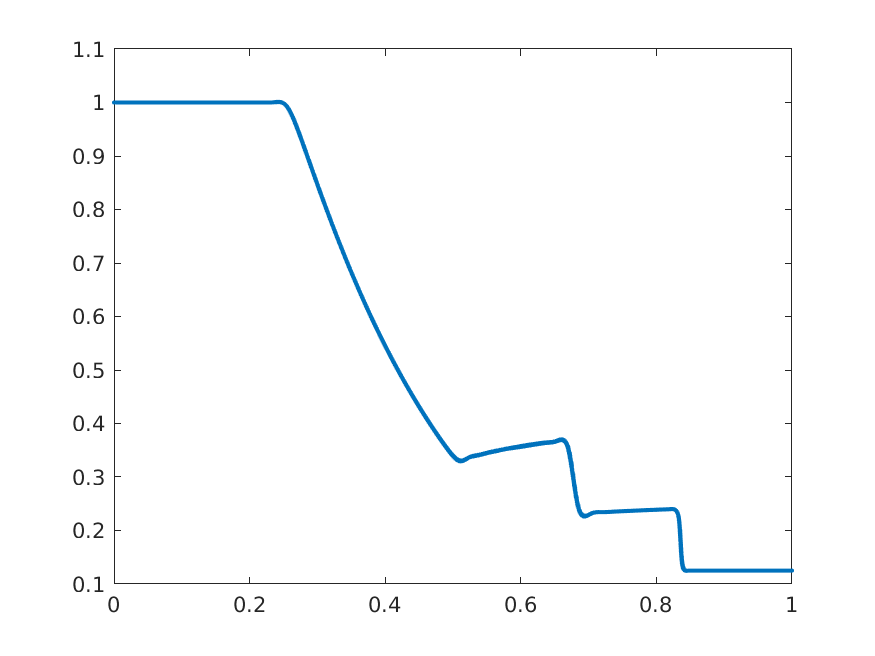}
  \includegraphics[width=0.32\linewidth]{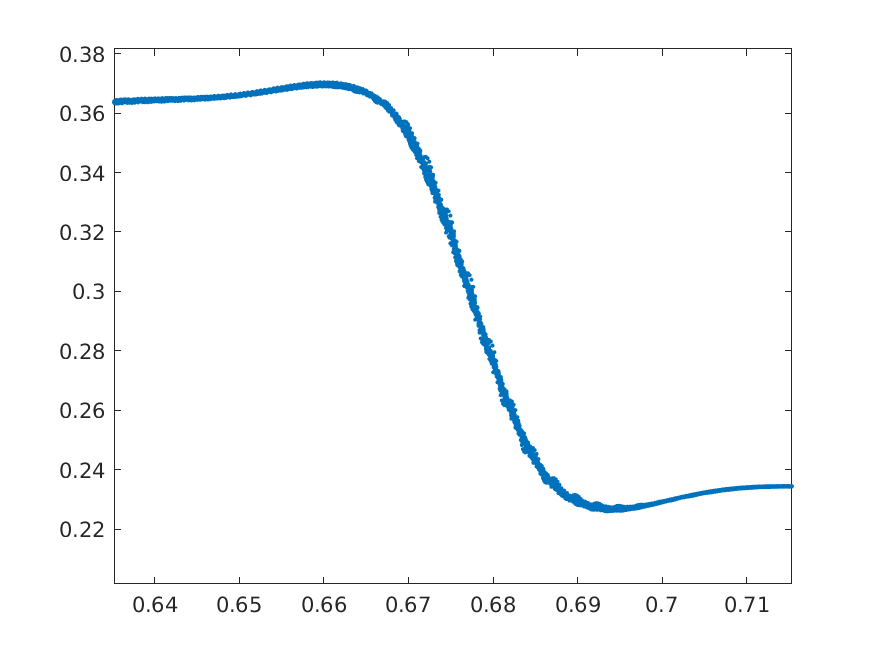}
  \includegraphics[width=0.32\linewidth]{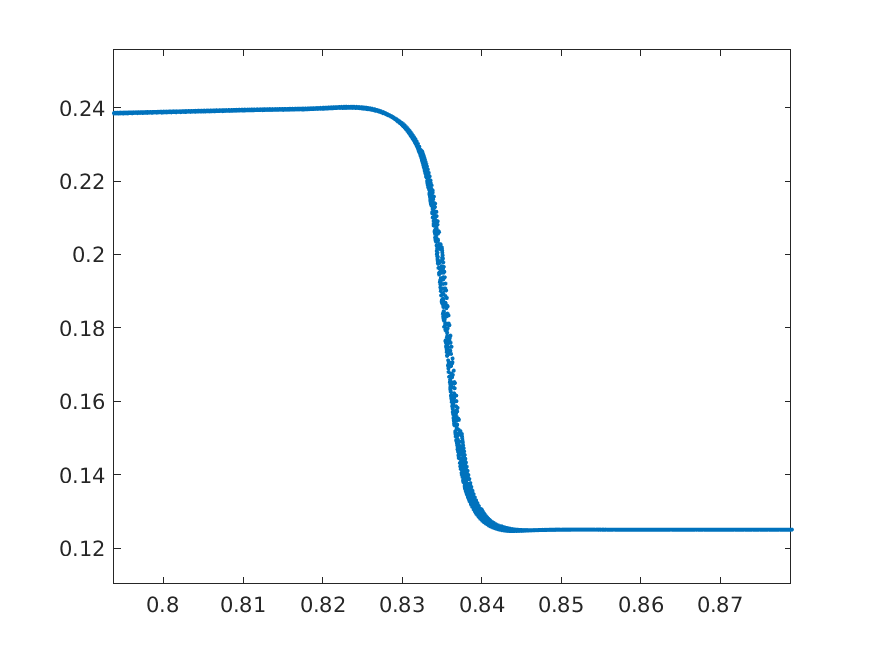}
\end{center}
	\caption{Radial Sod solutions at $t=0.2$ with the no-ghost \CWZb3 reconstruction.
		Density as a function of cell center from the origin for all cells.
		Whole solution (left), zoom on the contact (middle) and on the shock (right).}
\label{fig:sod:2dradial}
\end{figure}

In Fig.~\ref{fig:sod:2d} we compare the solutions at $t=0.2$ 
computed on a grid of $400\times400$ cells,
with and without using ghost cells. 
The solution is colored by pressure and 
25 equispaced contour lines of the density field, from $0.04$ to $1.0$,
 are also shown (grayscale colorbar),
so that the type of wave can be easily recognized.
All solutions were computed in the first quadrant only,
but the \CWENOZ3 one is shown reflected to the left to ease the comparison.

Almost no difference can be appreciated between the two solutions
and even the small artifacts that render the numerical solution non symmetrical appear identical in both schemes.
Further, in Fig.~\ref{fig:sod:2dradial},
we plot the density computed without using ghost cells
against the cell center distance from the origin.
An almost perfect radial symmetry is observed, despite the fact that 
the boundary cells are reconstructed with a different algorithm than the bulk ones.

\begin{figure}
\begin{center}
	\includegraphics[width=\linewidth]{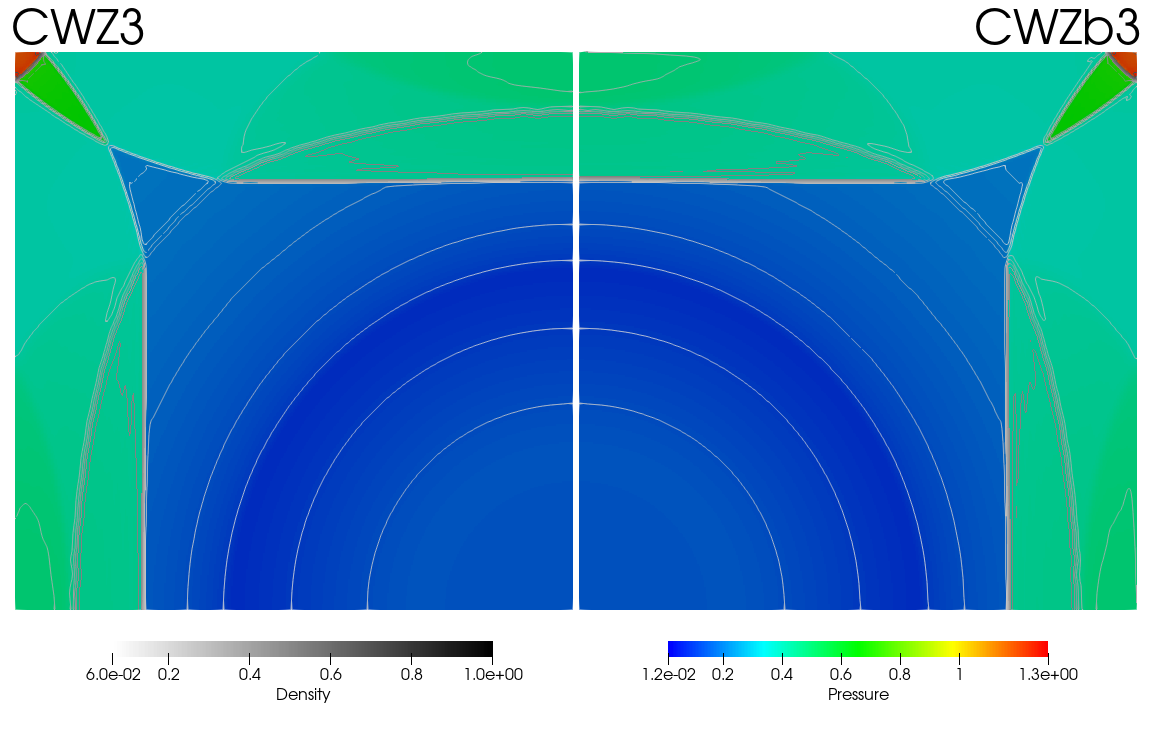}
\end{center}
\caption{Radial Sod solutions at $t=6$.}
\label{fig:sod:2dradial:long}
\end{figure}

The solution of this problem after $t=0.2$
sees the reflection of the cylindrical shock wave with the outer walls
and its later interaction with the expanding contact.
The reflected curved shock interacts with itself exactly at the upper-right corner 
at $t\simeq0.56$;
in Fig.~\ref{fig:sod:2dradial:long}, 
we show the solution at $t=0.6$, just after this event.
In this way we are testing the numerical schemes on reflecting a non planar shock wave on a wall.
Even more importantly, 
we are stressing the reconstruction procedure in the corner cell, 
since the shock convergence happening there 
implies that for some timesteps 
there is no non-discontinuous stencil available to the reconstruction procedure for the corner cell.
Here too, 
no appreciable difference is visible between the solutions 
computed with the two reconstruction schemes, 
showing that not using ghost cells in the reconstruction does not impair the numerical scheme.

\subsection{Implosion problem}
Next we consider the problem of a diamond-shaped converging shock
proposed in \cite{HuiLiLi:implosion}.
It is computed on a quarter plane with symmetry boundary conditions,
so that an oblique shock interacts with the boundary for a long time
in the initial stages of the evolution
and later the converging shock hits exactly the origin and is reflected back from there,
testing again the non-oscillatory properties of the reconstruction in a corner cell
when a strong wave impinges there.

The test is set in the square domain $[0,0.3]^2$
with reflective boundary conditions on all four sides: 
those at $x=0$ and $y=0$ represent symmetry lines, while the other two are physical solid walls.
The initial condition has zero velocity everywhere and $\rho=1,p=1$ in the outer region ($x+y>0.15$)
and $\rho=0.125,p=0.14$ in the interior one.
A useful reference for this test is \cite{LW:03} and the first author's website cited therein.
We show the solution computed with a grid of $800\times800$ cells;
the final time was set to $t=2.5$, saving snapshots
every $0.005$ until $t=0.1$ and every $0.1$ afterwards.

\begin{figure}
\begin{center}
\includegraphics[width=\linewidth]{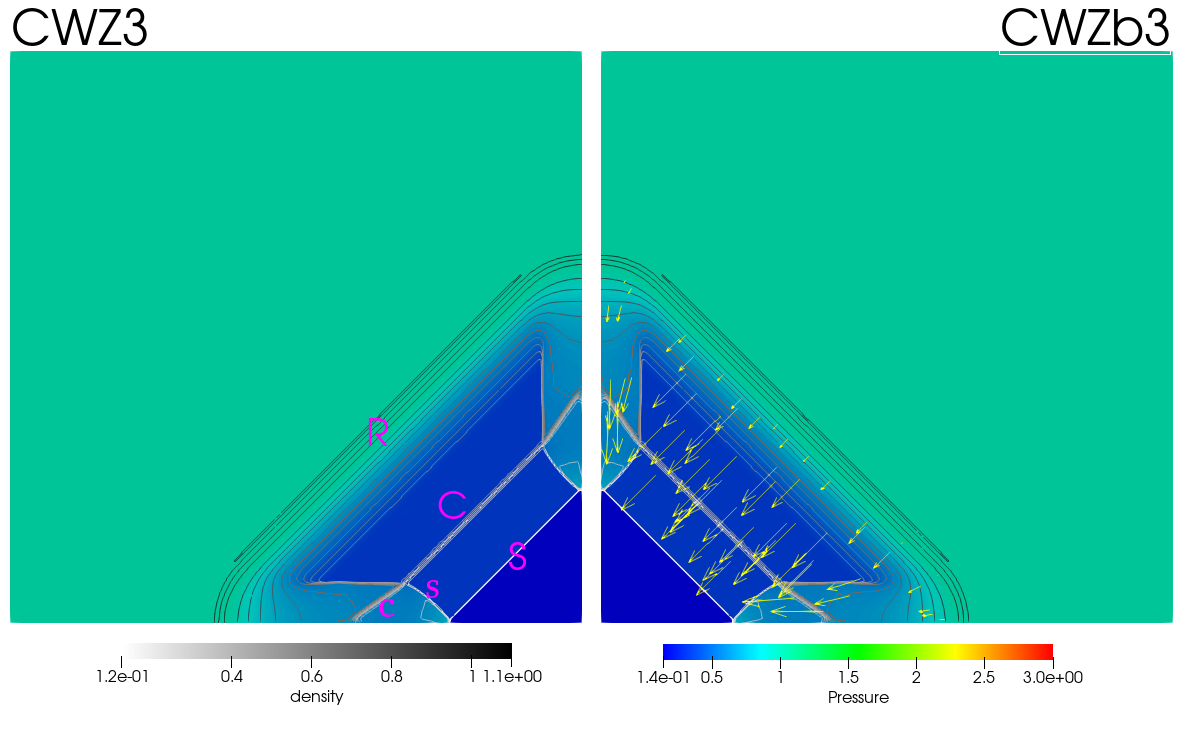}
\end{center}
\caption{Implosion test at $t=0.03$.
The rainbow colorbar is for pressure, 
the grayscale one is for the density isolines.
In the right panel, the arrows represent the velocity.
In the left panel,
the main shock (S), contact (C) and rarefaction (R) 
are indicated with capital letters,
some secondary waves with small letters.
}
\label{fig:impl:init}
\end{figure}

Fig.~\ref{fig:impl:init} shows both solutions 
in an early stage of the evolution,
at $t=0.03$.
Here and in all subsequent figures, we have reflected to the left the solution computed with ghosts.
In the early stages of the evolution
the initial discontinuity gives rise to 
a shock (indicated with ``S'' in the left panel) 
and a contact (``C''),
both moving towards the origin,
and to a rarefaction (``R'') that moves outwards.
At the boundary, the shock is reflected
and the reflected waves interact with the incoming contact
(``s'' and ``c'' in the figure).
In the right panel, the gas velocity is represented with arrows;
notice the fast wind directed towards the origin
blowing along the coordinate axis.

\begin{figure}
	\begin{center}
		\includegraphics[width=\linewidth]{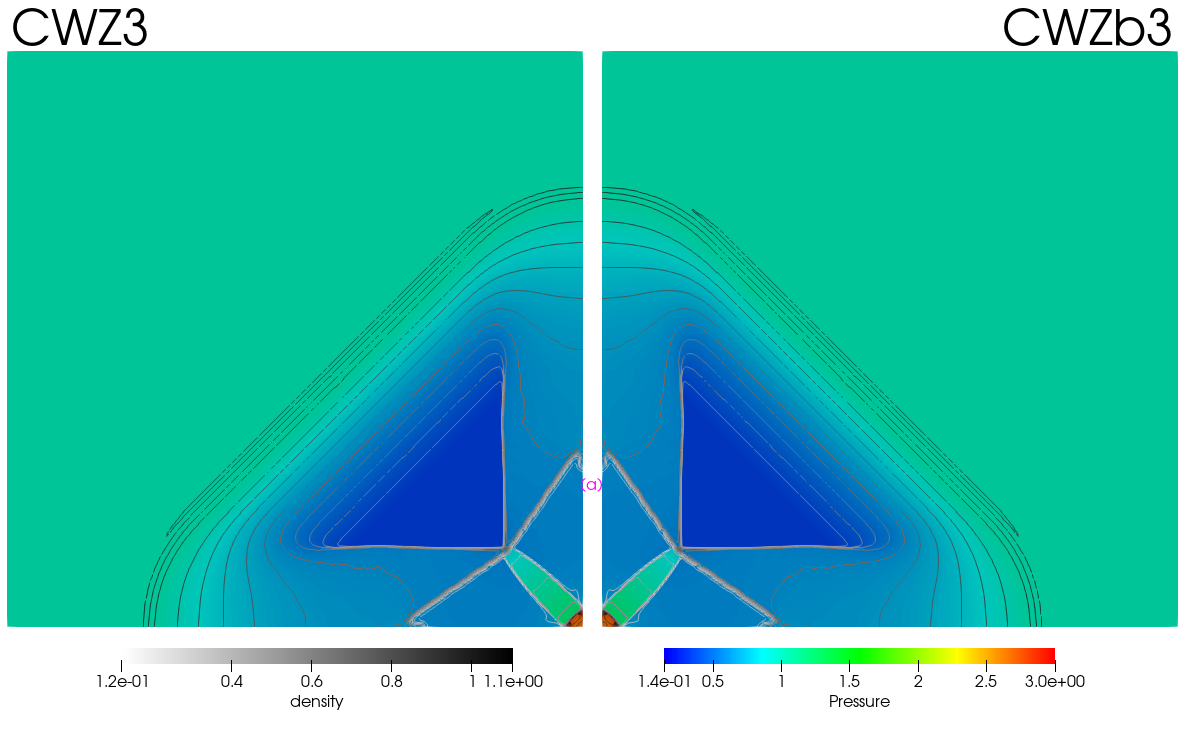}
	\end{center}
	\caption{Implosion test at $t=0.06$.}
	\label{fig:impl:shockcorner}
\end{figure}

Later the main shock and the reflected shocks 
converge in the origin,
hit there head to head and are bounced back outwards.
The snapshot reported in Fig.~\ref{fig:impl:shockcorner} 
is taken at $t=0.06$, just after this event.
Here it is important to observe that no spurious waves
and no difference among the two schemes can be observed close to the origin,
testifying that the reconstruction procedure in the corner cells is able to employ correctly the constant polynomial when the waves are very close to the corner.

The reflected contact, 
further deformed by the interaction with the expanding reflected shock, 
is being deformed by the 
wind blowing along the coordinate axes;
this is quite visible now at the point indicated by (a) in the figure.
Also this feature of the flow is computed symmetrically by both the ghosted and the no-ghost schemes.

\begin{figure}
	\begin{center}
		\includegraphics[width=\linewidth]{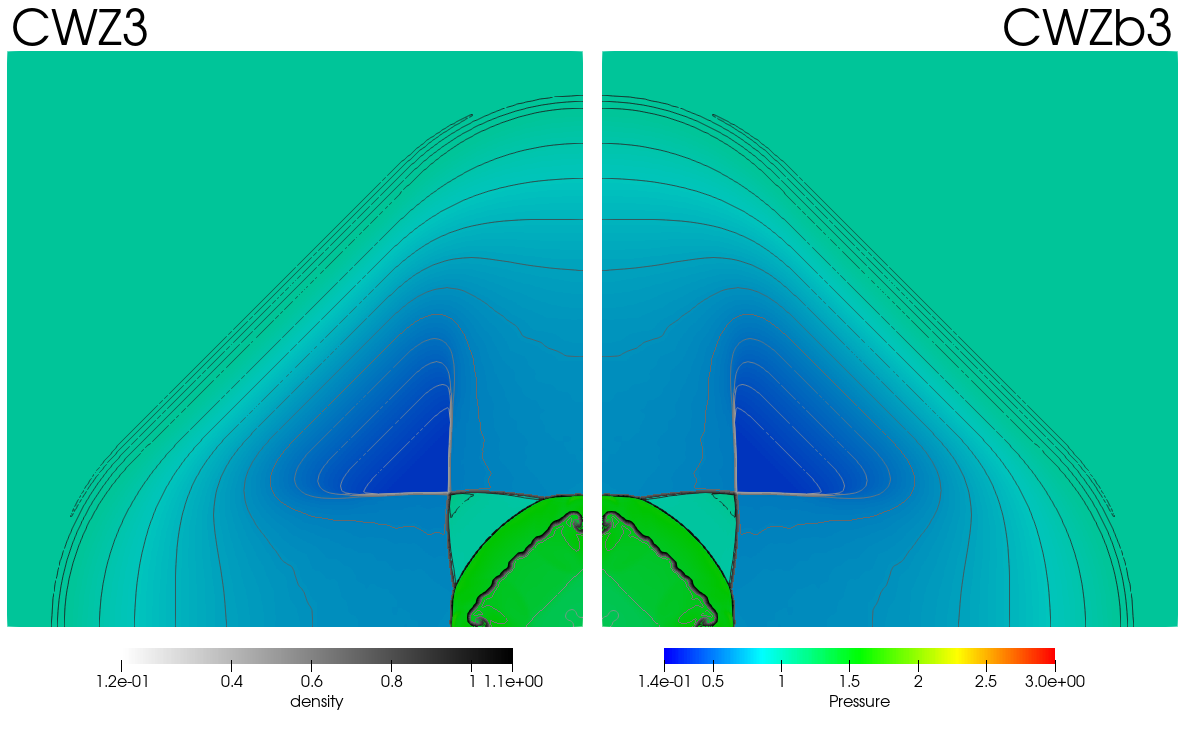}
	\end{center}
	\caption{Implosion test at $t=0.1$.}
	\label{fig:impl:01}
\end{figure}

In Fig.~\ref{fig:impl:01} we show the solution at time $t=0.1$.
At this time the rarefaction is still moving outwards, 
while the rounded shock bounced back from the origin 
has overcome the incoming contact,
which shows its physical instability coming from the deformations along the coordinate axis.

\begin{figure}
\begin{center}
	\includegraphics[width=\linewidth]{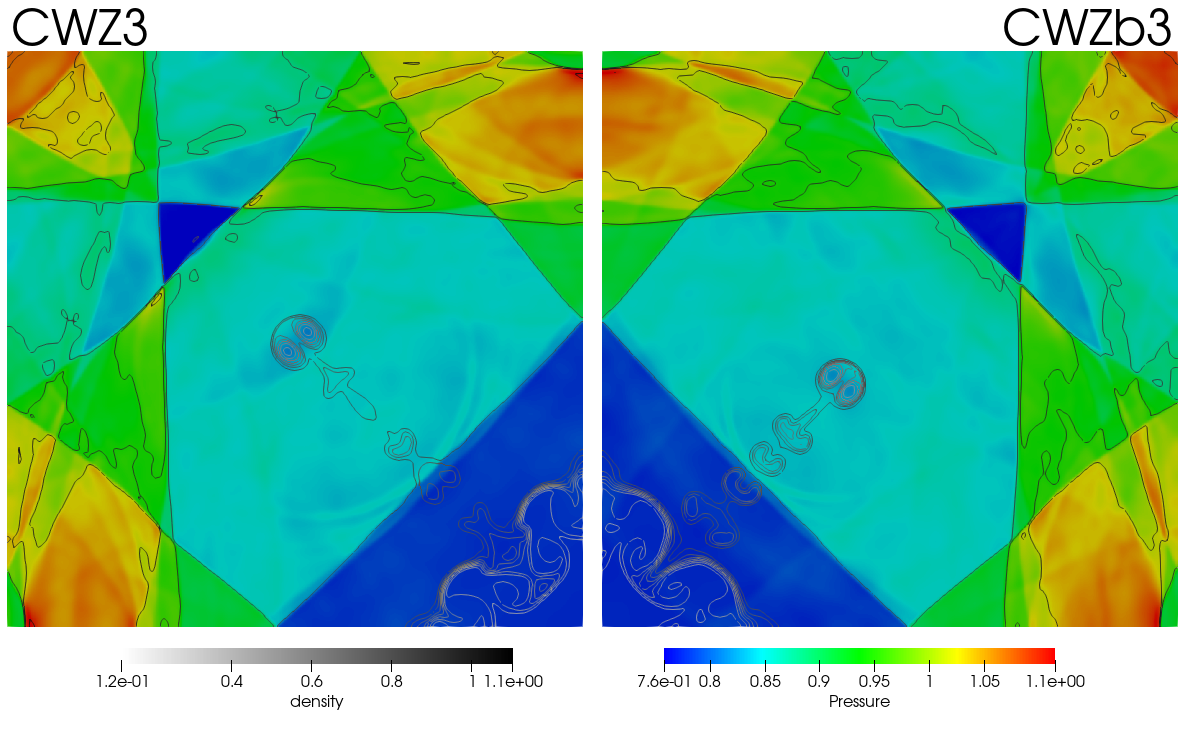}
\end{center}
\caption{Implosion test $t=2.5$. Note that this figure has a different colorbar than the previous ones.}
\label{fig:impl:25}
\end{figure}

When the evolution is computed for long times, there is no consensus among the different schemes
about the form of the bubbles near the origin and along the main diagonal \cite{LW:03}.
In Fig.~\ref{fig:impl:25} we report the solutions computed for $t=2.5$ 
(note the different pressure colorbar than in the previous ones).
At this time, many waves reflections, refractions and interactions have taken place
and the solution exhibit a quite complex pattern.
The two schemes compute the main waves and the pressure field almost identically,
but differ from each other in the bubbles, 
which appear to be physically unstable slip lines 
and therefore it is quite natural that different schemes can represent them in a different way.

\subsection{Shock-bubble interaction}
The last computation that we show is the shock-bubble interaction problem from \cite{CadaTorrilhon09}.
Here a right-moving shock 
hits a standing bubble of gas at low pressure.
In the computational domain $[-0.1, 1.6] \times [-0.5, 0.5]$, 
three distinct areas are considered: 
the post-shock region (A) for $x < 0$, 
the bubble (B) of center $(0.3, 0.0)$ and radius $0.2$ 
and the pre-shock region (C) of all points with $x > 0$ and not in (B). 
The initial data are
$\rho = \frac{11}{3}, u = 2.7136021011998722, p = 10.0$ in A,
$\rho = 0.1, u = 0.0, p = 1.0$ in B
and
$\rho = 1.0, u = 0.0, p = 1.0$ in C.
The vertical velocity $v$ is set to zero everywhere.
Boundary conditions are of Dirichlet type on the left (equal to the initial data), 
free-flow on the right, solid walls on $y = \pm 0.5$. 
The symmetry in the $y$ variable permits a half
domain computation (with $y > 0$) considering symmetry boundary conditions at $y = 0$.

\begin{figure}
\begin{center}
	\includegraphics[width=\linewidth]{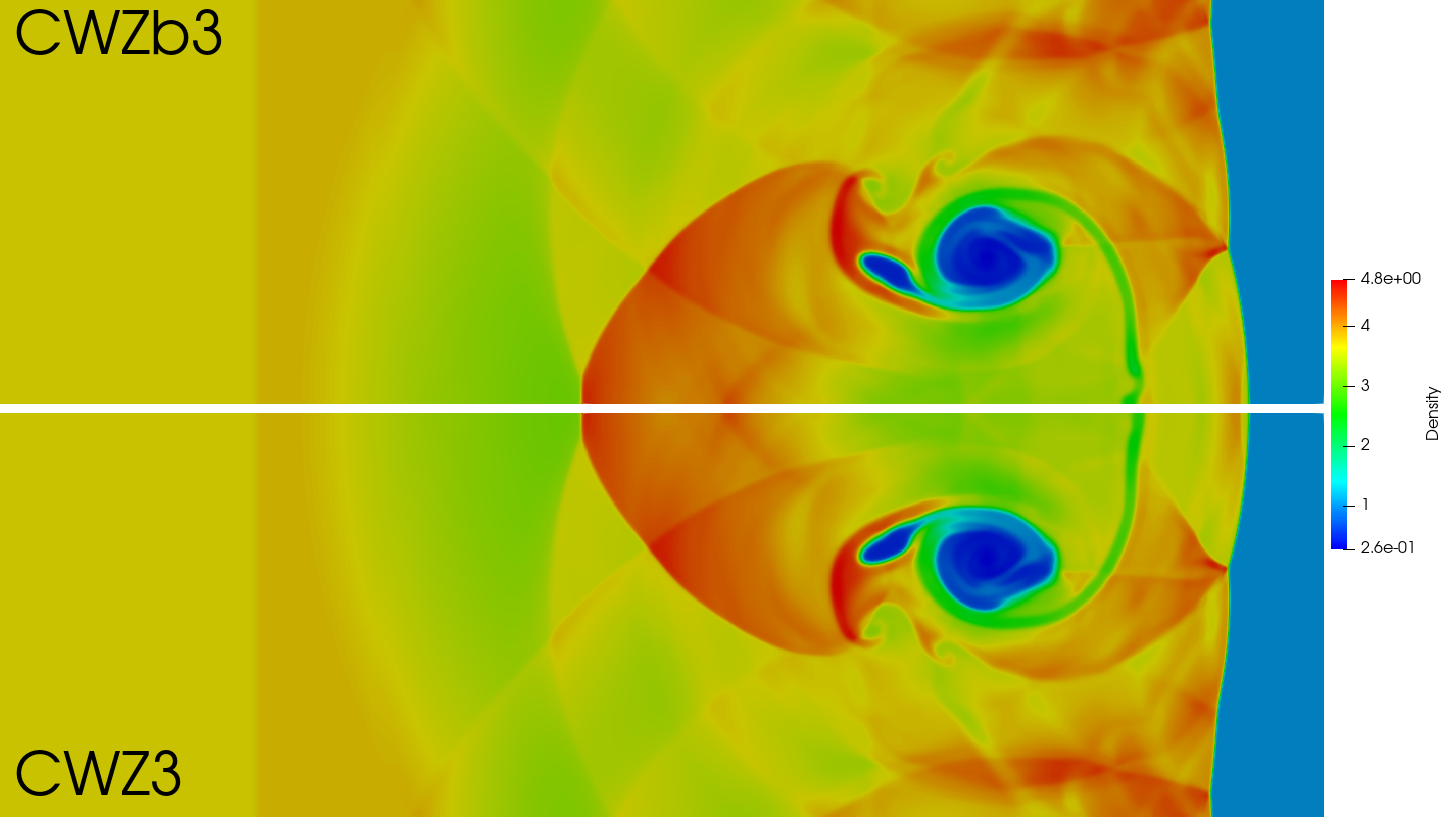}
	\\[1mm]
	\includegraphics[width=\linewidth]{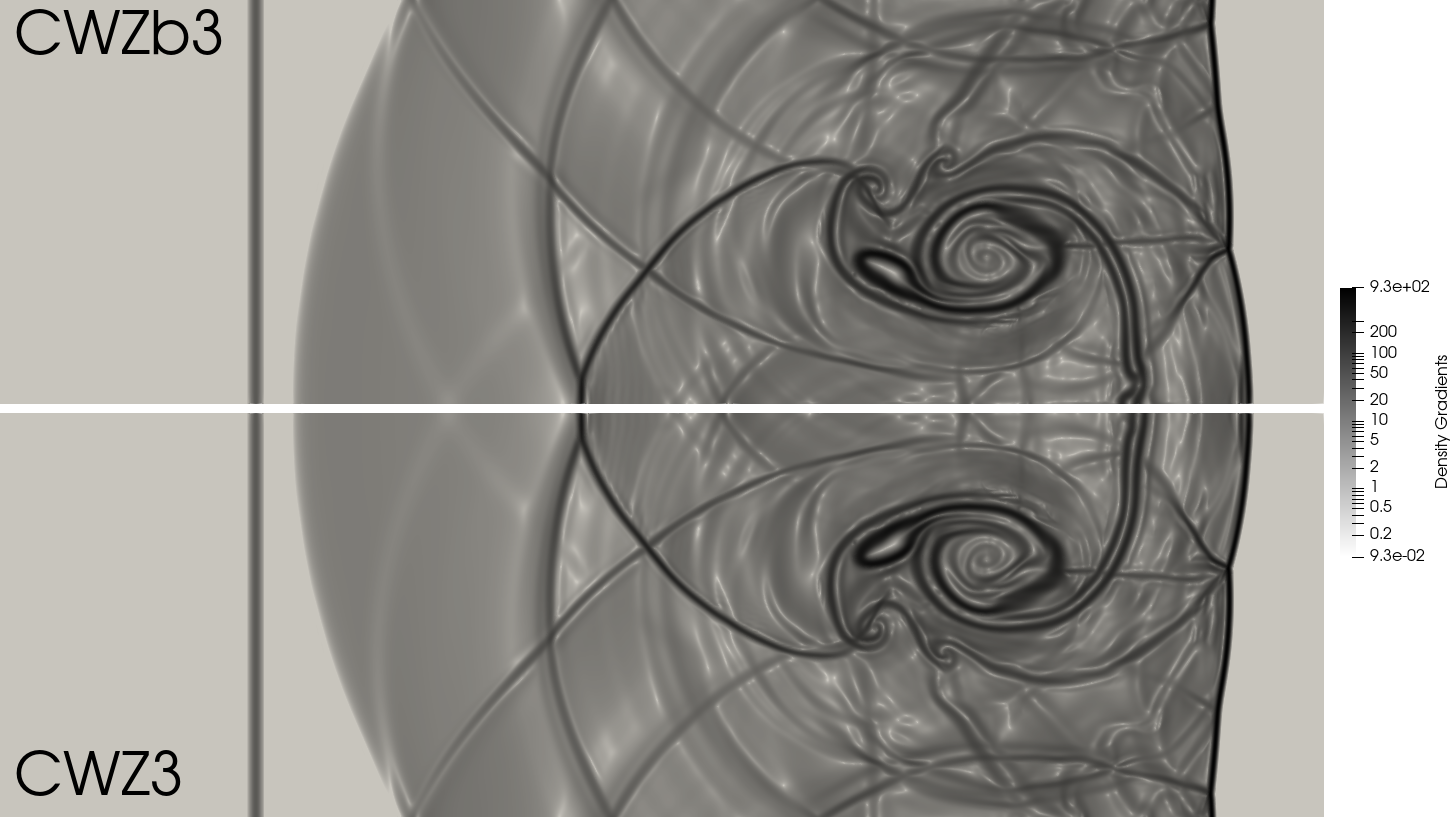}
\end{center}
\caption{Shock-bubble test: at $t=0.4$.
The top panel is colored by density.
The bottom panel is a numerical Schlieren plot 
obtained by representing the magnitude of the density gradient
in a grayscale logarithmic colorbar.
}
\label{fig:sb}
\end{figure}

The shock, in its movement towards the right, sets in motion, compresses and deforms the bubble; 
it interacts with it and the resulting refracted shocks are bounced back towards $y=0$ by the outer walls,
giving rise to a very complex interaction pattern.
The bubble is an unstable pattern 
and computing the solution with different schemes 
or different grid resolution will give it different final shapes.
In Fig.~\ref{fig:sb} we show the solutions at the final time of the computation, $t=0.4$;
the \CWENOZ3 solution is reflected along the symmetry axis in order to ease the comparison
with the no-ghost solution computed with \CWZb3. 
Here the main difference can be seen in the portion of the bubble that remains attached to the symmetry plane. All the other waves, including the central portion of the bubble, are almost identical in both computations.

\section{Conclusions and perspectives}
\label{sec:concl}

In this paper we have pursued further the approach of \cite{KolbSemplice:boundary}
for reconstructions in finite volume schemes without using ghost cells.
While the main motivation there was the application in internal nodes of networks,
where it is difficult to prescribe appropriate extrapolations,
here we have focused on the accuracy of the reconstructions and 
on the extension to higher dimensions.
All the proposed reconstructions in the boundary cells are based 
on reconstruction stencils that are not symmetric 
(but extended only inwards)
and nevertheless the optimal accuracy on smooth flows can be achieved.

Regarding the first point, we have proposed the employment of Z-weights
instead of the Jiang-Shu nonlinear weights,
obtaining a reconstruction that has less numerical diffusion
than the one of \cite{KolbSemplice:boundary}.
Perhaps more importantly,
the novel reconstruction can reach the optimal third order of convergence
in a more ample subset of the parameter space
and in particular can employ a very small $\epsilon$
without needing the infinitesimal weight of the constant polynomial
to be $\Ogrande(\DX^2)$.

In the second part of the paper,
we have proposed a (non dimensionally-split) reconstruction 
for two-dimensional Cartesian grids.
This has been compared with the \CWENOZ\ approach of \cite{CSV19:cwenoz}
that uses ghost cells,
showing that the employment of ghost cells can be entirely avoided
without affecting the quality of the computed solutions.

The approach without using ghosts appears thus to be quite promising
since setting them is always a tricky point in numerical schemes;
it seems interesting to pursue further this line of research 
towards higher accuracy or for moving boundaries like in piston problems.

\begin{acknowledgements}
The authors wish to thank Dr. Giuseppe Visconti for useful discussions
during the development of this paper.

This work was funded by Ministero dell'Università dello Stato Italiano,
under the PRIN project 2017KKJP4X
``Innovative numerical methods for evolutionary partial differential equations and applications''.
\end{acknowledgements}

\paragraph{Conflict of interest} 
Matteo Semplice, 
Elena Travaglia and
Gabriella Puppo
declare that they have no conflict of interest.

\bibliographystyle{plain}
\bibliography{CWBpaper.bib}

\end{document}